\documentclass[12pt]{article}
\usepackage[]{amsmath,amssymb}
\usepackage{amscd}
\usepackage{latexsym}
\usepackage{cite}

\newtheorem{definition}{Definition}[section]
\newtheorem{theorem}[definition]{Theorem}
\newtheorem{lemma}[definition]{Lemma}
\newtheorem{corollary}[definition]{Corollary}

\newtheorem{example}[definition]{Example}

\newtheorem{problem}[definition]{Problem}
\newtheorem{note}[definition]{Note}

\newtheorem{proposition}[definition]{Proposition}

\def\N{\mathbb N}

\def\Z{\mathbb Z}
\def\F{\mathbb F}

\begin{document}
\title{\bf  
The algebra $U_q({\mathfrak{sl}_2})$ in disguise
}
\author{
Sarah Bockting-Conrad
and
Paul Terwilliger
}
\date{}

\maketitle
\begin{abstract}
We discuss a connection between 
the algebra
$U_q({\mathfrak{sl}_2})$ and the tridiagonal pairs
of $q$-Racah type. To describe the connection,
let $x,y^{\pm 1},z$ denote the
equitable generators for 
$U_q({\mathfrak{sl}_2})$. Let $U^\vee_q$ denote
the subalgebra of 
$U_q({\mathfrak{sl}_2})$ generated by
$x,y^{-1},z$. Using
a tridiagonal pair of $q$-Racah type we construct
two finite-dimensional $U^\vee_q$-modules.
The constructions yield two nonstandard presentations of
 $U^\vee_q$ by generators and relations.
These presentations are investigated in detail.


 \end{abstract}

\section{Introduction}

\noindent This paper is about
a connection
between the algebra
 $U_q({\mathfrak{sl}_2})$ 
and
a linear-algebraic object called 
a tridiagonal pair \cite{TD00}.
To describe the connection,
we briefly recall the algebra 
$U_q({\mathfrak{sl}_2})$ 
\cite{jantzen,kassel}.
We will use the 
equitable presentation,
which was introduced in
\cite{equit}.
Given a field $\F$,
the $\F$-algebra 
$U_q({\mathfrak{sl}_2})$
 has generators
$x,y^{\pm 1},z$ and relations $yy^{-1}=y^{-1}y=1$,
\begin{eqnarray*}
\frac{qxy-q^{-1}yx}{q-q^{-1}} = 1,
\qquad
\frac{qyz-q^{-1}zy}{q-q^{-1}} = 1,
\qquad
\frac{qzx-q^{-1}xz}{q-q^{-1}} = 1.
\end{eqnarray*}
Let $U^\vee_q$ denote the subalgebra of
$U_q({\mathfrak{sl}_2})$ generated by
$x,y^{-1},z$. 
We now briefly turn our attention to tridiagonal pairs.
Let $V$ denote a vector space over $\F$ 
with finite positive dimension.
Roughly speaking, a tridiagonal pair on $V$
is a pair of diagonalizable $\F$-linear maps on $V$,
each acting on the eigenspaces of the other one
in a block-tridiagonal fashion.
There is a general family of tridiagonal pairs
said to have $q$-Racah type
\cite[Definition~3.1]{IT:qRacah}.
Let ${\bf A}, {\bf A^*}$
denote a tridiagonal pair
on $V$ that has 
$q$-Racah type.
Associated with this tridiagonal pair are two
additional 
$\F$-linear maps
$K:V\to V$
\cite[Definition~3.1]{bockting2}
and $B:V\to V$
\cite[Definition~3.2]{bockting2}, which are roughly
described as follows.
The maps $K,B$ are diagonalizable, and their
eigenspace decompositions are among
the split decompositions
of $V$ with respect to ${\bf A}, {\bf A^*}$ 
\cite[Section~4]{TD00}.
 According to
\cite[Lemma~3.6, Theorem~9.9]{bockting2}, there exists
$a \in \F \backslash \lbrace 0,1,-1\rbrace$ such that
\begin{eqnarray*}
&&
\frac{qK {\bf A} - q^{-1} {\bf A} K}{q-q^{-1}} = a K^2 + a^{-1} I,
\qquad \qquad 
\frac{qB {\bf A} - q^{-1} {\bf A} B}{q-q^{-1}} = a^{-1} B^2 + a I,
\\
&&
\qquad \quad
aK^2 - \frac{a^{-1} q-aq^{-1}}{q-q^{-1}} KB 
- \frac{a q-a^{-1} q^{-1}}{q-q^{-1}} BK  
+
a^{-1} B^2 = 0.
\end{eqnarray*}
The main purpose of this paper is to explain 
what the above three equations have to do with
 $U_q({\mathfrak{sl}_2})$. 
To summarize the answer,
consider an $\F$-algebra defined by
generators and relations; the generators
are symbols $K,B,{\bf A}$ and the relations
are the above three equations.
We are going to show that this algebra is isomorphic to
 $U^\vee_q$.
The existence of this isomorphism
yields a presentation of
$U^\vee_q$ by generators and relations. 
There is another presentation of 
$U^\vee_q$ with similar features, obtained by
inverting $q,a,K,B$.
 In our main results,
we describe these two presentations in a comprehensive way
that places them
 in a wider context.
  We will summarize
 our results shortly.

\medskip
\noindent 
In the theory of 
$U_q({\mathfrak{sl}_2})$, it is convenient
to introduce elements $\nu_x, \nu_y,\nu_z$ 
defined by
\begin{eqnarray*}
\nu_x = q(1-yz),
\qquad 
\nu_y = q(1-zx),
\qquad
\nu_z = q(1-xy).
\end{eqnarray*}
One significance of 
$\nu_x, \nu_y,\nu_z$ is that
\begin{eqnarray*}
&& x \nu_y = q^2 \nu_y x, \qquad \qquad x \nu_z = q^{-2} \nu_z x,
\\
&& y \nu_z = q^2 \nu_z y, \qquad \qquad y \nu_x = q^{-2} \nu_x y,
\\
&& z \nu_x = q^2 \nu_x z, \qquad \qquad z \nu_y = q^{-2} \nu_y z.
\end{eqnarray*}
\noindent We now summarize our results.
Fix a nonzero $a \in \F$ such that $a^2\not=1$.
Define elements $X,Z$ of
$U_q({\mathfrak{sl}_2})$ by
\begin{eqnarray*}
&&X = a^{-2}x + (1-a^{-2})y^{-1},
\qquad \qquad 
Z = a^{2}z + (1-a^{2})y^{-1}.
\end{eqnarray*}
One readily verifies that 
\begin{eqnarray}
\label{eq:match}
a^{-1}x+az=
aX+a^{-1}Z.
\end{eqnarray}
Let $A$ denote the common value of
(\ref{eq:match}).
We consider two antiautomorphisms
of $U_q({\mathfrak{sl}_2})$, denoted
$\tau$ and $\dagger$.
The map $\tau$ 
swaps $x\leftrightarrow z$
and fixes
$y^{\pm 1}$.
The map $\dagger$ 
swaps $x\leftrightarrow Z$,
  $z\leftrightarrow X$ and
fixes each of $A,y^{\pm 1}$.
The composition of
$\dagger$ and $\tau$ is 
an automorphism 
of $U_q({\mathfrak{sl}_2})$ which we
denote by
$\sigma$.
The map $\sigma$ sends $x\mapsto X$, $z\mapsto Z$ and
fixes $y^{\pm 1}$.
We give four presentations of 
$U_q({\mathfrak{sl}_2})$  by generators and relations.
They are 
 the Chevalley presentation, with generators
$e,f,k^{\pm 1}$;
 the equitable presentation, with generators 
$x,y^{\pm 1},z$;
 the modified Chevalley presentation, with
generators $\nu_x, y^{\pm 1}, \nu_z$;
 and the invariant presentation, with generators
$A,y^{\pm 1}$.
We show how these presentations are related.
We show how the maps
$\tau,\dagger,\sigma$ look 
in each presentation.
We show how the Casimir element looks in each presentation.
We then consider the subalgebra $U^\vee_q$.
We give four presentations of
$U^\vee_q$ by generators and relations.
For these presentations 
the generators are
(i) $x,y^{-1},z$;
(ii) $X,y^{-1},Z$;
(iii) $x,X,A$;
(iv) $z,Z,A$.
We explain how these presentations are related to
each other and the four presentations of 
$U_q({\mathfrak{sl}_2})$.
We show that $U^\vee_q$ is invariant under $\tau,\dagger,\sigma$.
We show how 
$\tau,\dagger,\sigma$
act on the generating sets
(i)--(iv).
We compare the finite-dimensional modules for
$U^{\vee}_q$ and
$U_q({\mathfrak{sl}_2})$.
Given a tridiagonal pair of $q$-Racah type we
construct two finite-dimensional 
$U^{\vee}_q$-modules; in these constructions
the presentations (iii), (iv)
come up 
in a natural way. Presentation (iv) 
is the one we discussed at the outset, if
we identify
$z\mapsto K$, $Z\mapsto B$, $A\mapsto {\bf A}$.
Given a tridiagonal pair of $q$-Racah type,
there is an associated linear map  $\psi$ called the
 Bockting operator
\cite{twocom,bockting2}.
The map $\psi$ 
plays the following role in our
theory.
Consider our two $U^\vee_q$-modules constructed from
the given tridiagonal pair. We show that each of these 
$U^\vee_q$-modules extends to a 
$U_q({\mathfrak{sl}_2})$-module.
We show that on one of the 
$U_q({\mathfrak{sl}_2})$-modules
$\nu_x$ acts as
$a^{-1}\psi$, and on the other 
$U_q({\mathfrak{sl}_2})$-module
 $\nu_z$ acts as
$a\psi$.
\medskip

\noindent
At the end of the paper we list some open problems concerning
$U^\vee_q$ and related topics.

\section{
The algebra $U_q({\mathfrak{sl}_2})$ }

\noindent 
\noindent We now begin our formal argument.
We adopt the following conventions.
An algebra is meant to be associative and have a 1.
A subalgebra has the same 1 as the parent algebra.
Let $\F$ denote a field and fix a nonzero $q \in \F$
such that $q^4 \not= 1$. 
\medskip

\noindent We now recall the quantum  enveloping algebra
$U_q({\mathfrak{sl}_2})$. For background information on this topic
we refer the reader to
the books by Jantzen
\cite{jantzen}
and
Kassel
\cite{kassel}.

\begin{definition}
\label{def:uqchev}
\rm 
\cite[Section~1.1]{jantzen}.
Let $U_q({\mathfrak{sl}_2})$ denote the $\F$-algebra defined by
generators $e$, $f$, $k^{\pm 1}$ and relations
\begin{eqnarray}
&&\qquad \qquad kk^{-1} =  k^{-1}k = 1,
\label{eq:k}
\\
&&kek^{-1} = q^2 e, \qquad \qquad kfk^{-1} = q^{-2}f,
\label{eq:ke}
\\
&&\qquad \quad \;\;  ef - fe = \frac{k-k^{-1}}{q-q^{-1}}.
\label{eq:ef}
\end{eqnarray}
\end{definition}

\noindent
The presentation 
of $U_q({\mathfrak{sl}_2})$ given in
Definition
\ref{def:uqchev}
is called the {\it Chevalley presentation}.

\medskip

\noindent Recall the natural numbers
$\N = \lbrace 0,1,2,\ldots \rbrace$ and
integers $\Z=\lbrace 0,\pm 1,\pm 2,\ldots\rbrace$. 

\begin{lemma} 
{\rm \cite[Theorem~1.5]{jantzen}.} The following is a basis for the $\F$-vector
space 
$U_q({\mathfrak{sl}_2})$:
\begin{eqnarray}
\label{eq:uqbasis}
e^r k^s f^t \qquad \qquad r,t\in \N, \qquad s \in \Z.
\end{eqnarray}
\end{lemma}


\begin{definition}
\label{def:cas} 
\rm
{\rm \cite[Section~2.7]{jantzen}.}
Define
\begin{eqnarray*}
\Lambda &=& (q-q^{-1})^2 ef + q^{-1}k + q k^{-1}.
\end{eqnarray*}
$\Lambda$ is called the {\it Casimir element} of
$U_q({\mathfrak{sl}_2})$.
\end{definition}

\begin{lemma} 
\label{lem:centerbasis}
{\rm \cite[Lemma~2.7, Proposition~2.18]{jantzen}.}
The center of 
$U_q({\mathfrak{sl}_2})$ contains $\Lambda$.
This center has a basis
$\lbrace \Lambda^i\rbrace_{i \in \N}$,
provided that $q$ is not a root of unity.
\end{lemma}

\begin{lemma}
\label{lem:sigmasend}
For all nonzero $a \in \F$ there exists
a unique automorphism $\sigma_a$ of 
$U_q({\mathfrak{sl}_2})$ that sends
\begin{eqnarray}
e \mapsto a^{-2} e, \qquad \qquad f \mapsto a^2 f,
\qquad \qquad k^{\pm 1} \mapsto k^{\pm 1}.
\label{eq:sigmasend}
\end{eqnarray}
Moreover $\sigma^{-1}_a = \sigma_{a^{-1}}$.
\end{lemma}
\noindent {\it Proof:}  Use 
(\ref{eq:k})--(\ref{eq:ef}).
\hfill $\Box$ \\

\begin{lemma} For all nonzero $a \in \F$ the
automorphism $\sigma_a$ fixes $\Lambda$.
\end{lemma}
\noindent {\it Proof:} 
Use Definition
\ref{def:cas}  and 
(\ref{eq:sigmasend}).
\hfill $\Box$ \\

\noindent 
The automorphisms $\sigma_q$ and $\sigma_{q^{-1}}$ act on
$U_q({\mathfrak{sl}_2})$  as follows.

\begin{lemma} \label{lem:aqqi}
For $\xi \in 
U_q({\mathfrak{sl}_2})$, 
$\sigma_q$ sends $\xi \mapsto k^{-1}\xi k$
and $\sigma_{q^{-1}}$ sends
$\xi \mapsto k\xi k^{-1}$.
\end{lemma}
\noindent {\it Proof:} 
Use
(\ref{eq:ke}) and
(\ref{eq:sigmasend}).
\hfill $\Box$ \\

\section{The equitable presentation of 
$U_q({\mathfrak{sl}_2})$}

\noindent  In this section we recall the equitable presentation of
$U_q({\mathfrak{sl}_2})$ 
\cite{equit}. For background information on this topic
we refer the reader to
\cite{equit,tersym,uawe}. See also
\cite{alnajjar,bt,
uq}.

\begin{lemma} 
\label{def:uqequit}
{\rm
\cite[Theorem~2.1]{equit}.}
The $\F$-algebra 
$U_q({\mathfrak{sl}_2})$ has a presentation by generators
$x$, $y^{\pm 1}$,  $z$ and relations
\begin{eqnarray}
&& \qquad \qquad \qquad 
\qquad \qquad 
{yy^{-1} = y^{-1}y = 1},
\label{eq:yyi}
\\
&&
\frac{qxy-q^{-1}yx}{q-q^{-1}} = 1,
\qquad 
\frac{qyz-q^{-1}zy}{q-q^{-1}} = 1,
\qquad 
\frac{qzx-q^{-1}xz}{q-q^{-1}} = 1.
\label{eq:xyz}
\end{eqnarray}
An isomorphism with the presentation in Definition
\ref{def:uqchev} sends
\begin{eqnarray*}
x &\mapsto& k^{-1} - k^{-1} e q (q-q^{-1}),
\\
y^{\pm 1} &\mapsto& k^{\pm 1},
\\
z &\mapsto& k^{-1} +  f (q-q^{-1}).
\end{eqnarray*}
\noindent The inverse isomorphism sends
\begin{eqnarray*}
e &\mapsto & (1-yx)q^{-1}(q-q^{-1})^{-1},
\\
k^{\pm 1} &\mapsto& y^{\pm 1},
\\
f &\mapsto & (z-y^{-1})(q-q^{-1})^{-1}.
\end{eqnarray*}
\end{lemma}

\noindent
The presentation for
$U_q({\mathfrak{sl}_2})$
given in
Lemma 
\ref{def:uqequit}
is called the {\it equitable presentation}.

\begin{note}
\label{note:identify}
\rm From now on, we identify the
copy of
$U_q({\mathfrak{sl}_2})$ from Definition
\ref{def:uqchev} with the copy of
$U_q({\mathfrak{sl}_2})$ from Lemma
\ref{def:uqequit}, via the isomorphism in
Lemma
\ref{def:uqequit}.
\end{note}

\begin{lemma}
\label{lem:Uqbasis}
{\rm 
\cite[Lemma~10.7]{uawe}.
}
The following is a basis for the $\F$-vector
space
$U_q({\mathfrak{sl}_2})$:
\begin{eqnarray}
x^r y^s z^t \qquad \qquad r,t\in \N, \quad s \in \Z.
\label{eq:xyzbasis}
\end{eqnarray}
\end{lemma}

\noindent We now describe the Casimir element $\Lambda$ from the equitable
point of view.
\begin{lemma}
\label{lem:casequit}
{\rm \cite[Lemma~2.15]{uawe}.}
The element $\Lambda$ is equal to each of the following:
\begin{eqnarray*}
&&qx + q^{-1}y + qz-qxyz, \qquad \qquad q^{-1}x + qy + q^{-1}z-q^{-1}zyx,
\\
&&
qy + q^{-1}z + qx-qyzx, \qquad \qquad q^{-1}y + qz + q^{-1}x-q^{-1}xzy,
\\
&&
qz + q^{-1}x + qy-qzxy, \qquad \qquad q^{-1}z + qx + q^{-1}y-q^{-1}yxz.
\end{eqnarray*}
\end{lemma}

\noindent We now describe the automorphisms $\sigma_a$ of
$U_q({\mathfrak{sl}_2})$ from the equitable point of view.

\begin{lemma}
\label{lem:sigform}
For all nonzero $a \in \F$ the automorphism $\sigma_a$
sends 
\begin{eqnarray*}
x \mapsto a^{-2} x + (1-a^{-2})y^{-1},
\quad \qquad
y^{\pm 1} \mapsto y^{\pm 1},
\qquad \quad 
z \mapsto a^2 z + (1-a^2)y^{-1}.
\end{eqnarray*}
\end{lemma}
\noindent {\it Proof:} Use
(\ref{eq:sigmasend}) and the isomorphisms in
Lemma
\ref{def:uqequit}.
\hfill $\Box$ \\

\begin{lemma} 
\label{lem:yconj}
We have
\begin{eqnarray*}
&&y^{-1}xy = q^{-2}x + (1-q^{-2})y^{-1},
\qquad \qquad 
yxy^{-1} = q^{2}x + (1-q^{2})y^{-1},
\\
&&
yzy^{-1} = q^{-2}z + (1-q^{-2})y^{-1},
\qquad \qquad 
y^{-1}zy = q^{2}z + (1-q^{2})y^{-1}.
\end{eqnarray*}
\end{lemma}
\noindent {\it Proof:}  
These equations are reformulations of
the first two relations in
(\ref{eq:xyz}). They can also be obtained from
Lemma \ref{lem:aqqi}, by identifying $k=y$ and using
Lemma
\ref{lem:sigform}.
\hfill $\Box$ \\

\noindent We will be discussing antiautomorphisms, so let us
recall that concept.
Given an $\F$-algebra $\mathcal A$,
a map $\gamma:\mathcal A \to \mathcal A$ is called an
{\it antiautomorphism} whenever
$\gamma$ is an isomorphism of $\F$-vector spaces and
$(uv)^\gamma= v^\gamma u^\gamma$ for all $u,v \in \mathcal A$.

\begin{lemma} 
\label{lem:tau}
There exists a unique antiautomorphism $\tau$ of
$U_q({\mathfrak{sl}_2})$ that sends
\begin{eqnarray}
x \mapsto z, \qquad  \qquad 
y^{\pm 1} \mapsto y^{\pm 1}, \qquad 
\qquad z \mapsto x.
\label{eq:tausend}
\end{eqnarray}
Moreover $\tau^2=1$.
\end{lemma}
\noindent {\it Proof:}  
Let $S$ denote the
set of defining relations for
$U_q({\mathfrak{sl}_2})$ given in
Lemma
\ref{def:uqequit}. For each relation
$r \in S$ let $r'$ denote the
equation obtained by
inverting the order of
multiplication and swapping $x,z$.
The map $r\mapsto r'$ permutes $S$,
and therefore the antiautomorphism $\tau$ exists. 
The antiautomorphism $\tau$ is
unique since
$x,y^{\pm 1},z$ 
generate
$U_q({\mathfrak{sl}_2})$.
\hfill $\Box$ \\

\noindent We now consider how $\tau$ acts on $e,f,k^{\pm 1}$.

\begin{lemma} 
\label{lem:tauefk}
The antiautomorphism $\tau$ sends
\begin{eqnarray*}
e \mapsto -qkf, 
\qquad \qquad
f \mapsto -q^{-1}ek^{-1}, 
\qquad \qquad
k^{\pm 1} \mapsto k^{\pm 1}.
\end{eqnarray*}
\end{lemma}
\noindent {\it Proof:}  
Use
(\ref{eq:tausend}) and
the isomorphisms in
Lemma
\ref{def:uqequit}.
\hfill $\Box$ \\

\begin{lemma} The antiautomorphism $\tau$ fixes $\Lambda$.
\end{lemma}
\noindent {\it Proof:}  Use
Lemma
\ref{lem:casequit} and
(\ref{eq:tausend}).
\hfill $\Box$ \\

\begin{lemma}
\label{lem:tauprod}
For all nonzero $a \in \F$, 
$\tau \sigma_a = \sigma^{-1}_a \tau$.
\end{lemma} 
\noindent {\it Proof:}  
Use
Lemma
\ref{lem:sigmasend}
and Lemma
\ref{lem:tauefk}.
\hfill $\Box$ \\

\noindent We mention some identities for later use.

\begin{lemma}\label{lem:qsymrel}
We have
\begin{eqnarray}
&&
\frac{x^2y - (q^2+q^{-2})xyx+ yx^2}{(q-q^{-1})^2} = -x,
 \qquad
\frac{y^2x - (q^2+q^{-2})yxy+ xy^2}{(q-q^{-1})^2} = -y,
\label{eq:sym1}
\\
&&
\frac{y^2z - (q^2+q^{-2})yzy+ zy^2}{(q-q^{-1})^2} = -y,
\qquad
\frac{z^2y - (q^2+q^{-2})zyz+ yz^2}{(q-q^{-1})^2} = -z,
\label{eq:sym2}
\\
&&
\frac{z^2x - (q^2+q^{-2})zxz+ xz^2}{(q-q^{-1})^2} = -z,
\qquad
\frac{x^2z - (q^2+q^{-2})xzx+ zx^2}{(q-q^{-1})^2} = -x.
\label{eq:sym3}
\end{eqnarray}
\end{lemma}
\noindent {\it Proof:} We verify the equation on
the left in
(\ref{eq:sym1}).
Observe that
\begin{eqnarray*}
\frac{x^2y-(q^2+q^{-2})xyx + y x^2}{q-q^{-1}} = 
q^{-1}x \frac{qxy-q^{-1}yx}{q-q^{-1}}
- q\frac{qxy-q^{-1}yx}{q-q^{-1}}x.
\end{eqnarray*}
In the above equation, each fraction on the right is equal to 1,
so the whole expression on the right is equal to
$-(q-q^{-1}) x$. This yields the equation on the left
in
(\ref{eq:sym1}). The remaining equations are similarly verified.
\hfill $\Box$ \\

\section{A $\Z_2$-grading of
$U_q({\mathfrak{sl}_2})$}

\noindent We will be discussing the  group
$\Z_2 = \Z/2\Z$ of order 2.
Given an $\F$-algebra $\mathcal A$,
by a {\it $\Z_2$-grading} 
of $\mathcal A$ we mean a direct sum
decomposition
$\mathcal A=\sum_{i\in \Z_2} \mathcal A_i$ 
such that
$\mathcal A_i  \mathcal A_j \subseteq \mathcal A_{i+j}$ for $i,j\in \Z_2$.
In this case $\mathcal A_0$ is a subalgebra of $\mathcal A$.
\medskip

\noindent We now describe a $\Z_2$-grading of 
$U_q({\mathfrak{sl}_2})$.

\begin{definition}
\label{def:evenodd}
\rm
Referring to  the basis 
(\ref{eq:xyzbasis}) 
of
$U_q({\mathfrak{sl}_2})$, a basis element
$x^ry^sz^t$ will be called {\it even}
(resp. {\it odd}) whenever $r+s+t$ is even
(resp. odd). Let $U_{even}$ (resp. $U_{odd}$) denote
the subspace of 
$U_q({\mathfrak{sl}_2})$ spanned by
the even (resp. odd) elements in
(\ref{eq:xyzbasis}). 
By construction
\begin{eqnarray}
U_q({\mathfrak{sl}_2}) = U_{even}+U_{odd}
\qquad \qquad
{\mbox{\rm (direct sum).}}
\label{eq:grading}
\end{eqnarray}
\end{definition}

\begin{lemma} The decomposition
{\rm (\ref{eq:grading})} is a $\Z_2$-grading
of 
$U_q({\mathfrak{sl}_2})$.
In particular $U_{even}$ is a subalgebra of
$U_q({\mathfrak{sl}_2})$.
\end{lemma}
\noindent {\it Proof:} Recall the
defining relations
(\ref{eq:yyi}),
(\ref{eq:xyz}) for the equitable presentation of
$U_q({\mathfrak{sl}_2})$.
In these relations
every term has even degree.
\hfill $\Box$ \\

\begin{definition}
\label{def:uqprime}
\rm
Let $U_q'$ denote the subalgebra
of 
$U_q({\mathfrak{sl}_2})$ generated by $x, y, z$.
\end{definition}

\begin{lemma}
\label{lem:uprimebasis}
{\rm \cite[Lemma~10.8]{uawe}}.
 The following is a basis for the
$\F$-vector space $U_q'$:
\begin{eqnarray}
x^r y^s z^t \qquad \qquad r,s,t\in \N.
\label{eq:xyzbasispos}
\end{eqnarray}
\end{lemma}

\begin{definition}\label{def:evenoddprime}
\rm
Let $U'_{even}$ (resp. $U'_{odd}$) denote
the subspace of 
$U'_q$ spanned by
the even (resp. odd) elements in
(\ref{eq:xyzbasispos}). 
By construction
\begin{eqnarray}
U'_q = U'_{even}+U'_{odd}
\qquad \qquad
{\mbox{\rm (direct sum).}}
\label{eq:gradingpos}
\end{eqnarray}
\end{definition}

\begin{lemma} Referring to Definition
\ref{def:evenoddprime},
\begin{eqnarray}
U'_{even} = U'_q \cap U_{even},
\qquad \qquad 
U'_{odd} = U'_q \cap U_{odd}.
\label{eq:odprime}
\end{eqnarray}
Moreover the decomposition
{\rm (\ref{eq:gradingpos})} is a 
$\Z_2$-grading of $U_q'$.
\end{lemma}
\noindent {\it Proof:}  By construction.
\hfill $\Box$ \\

\begin{lemma}
For all nonzero $a \in \F$ the automorphism
$\sigma_a$ leaves invariant both
$U_{even}$ and $U_{odd}$.
However
$\sigma_a$ does not leave
$U'_q$ invariant, unless
 $a^2=1$.
\end{lemma}
\noindent {\it Proof:} 
Use Lemma \ref{lem:sigform}.
\hfill $\Box$ \\

\begin{lemma}
Each of the following is invariant under
the antiautomorphism $\tau$:
\begin{eqnarray*}
U_{even}, 
\qquad 
U_{odd}, 
\qquad 
U'_q,
\qquad 
U'_{even}, 
\qquad 
U'_{odd}.
\end{eqnarray*}
\end{lemma}
\noindent {\it Proof:} Use
Lemma \ref{lem:tau}.
\hfill $\Box$ \\

\begin{lemma}
\label{lem:casloc}
The Casimir element
$\Lambda$ is contained in  $U'_{odd}$.
\end{lemma}
\noindent {\it Proof:} Use
Lemma 
\ref{lem:casequit}.
\hfill $\Box$ \\

\begin{corollary} \label{cor:center}
Assume that $q$ is not a root of unity. Then
the center of 
$U_q({\mathfrak{sl}_2})$ is contained in
$U'_q$.
\end{corollary}
\noindent {\it Proof:} 
By Lemma
\ref{lem:centerbasis}
and
Lemma 
\ref{lem:casloc}.
\hfill $\Box$ \\

\section{The elements $\nu_x, \nu_y, \nu_z$
of $U_q({\mathfrak{sl}_2})$}

\noindent 
The elements
 $\nu_x, \nu_y, \nu_z$ of 
$U_q({\mathfrak{sl}_2})$ were introduced in
\cite{uawe}.
In this section, we use these elements
to obtain a generating set for the subalgebras
$U_{even}$ and $U'_{even}$.

\medskip

\noindent 
Reformulating the relations 
(\ref{eq:xyz}), we obtain
\begin{eqnarray*}
q(1-xy) = q^{-1}(1-yx),
 \qquad 
q(1-yz) = q^{-1}(1-zy),
\qquad 
q(1-zx) = q^{-1}(1-xz).
\end{eqnarray*}
Following 
\cite[Definition~3.1]{uawe}
we define
\begin{eqnarray}
&&\nu_x = q(1-yz)= q^{-1}(1-zy),
\label{eq:nux}
\\
&&\nu_y = q(1-zx)= q^{-1}(1-xz),
\label{eq:nuy}
\\
&&\nu_z = q(1-xy)= q^{-1}(1-yx).
\label{eq:nuz}
\end{eqnarray}
\noindent By Lemma
\ref{def:uqequit}
and Note
\ref{note:identify},
\begin{eqnarray}
\nu_x = -q(q-q^{-1})kf, \qquad \qquad \nu_z = (q-q^{-1})e.
\label{eq:xzef}
\end{eqnarray}

\begin{lemma} \label{lem:qcom}
{\rm \cite[Lemma~3.5]{uawe}. }
We have
\begin{eqnarray*}
&& x \nu_y = q^2 \nu_y x, \qquad \qquad x \nu_z = q^{-2} \nu_z x,
\\
&& y \nu_z = q^2 \nu_z y, \qquad \qquad y \nu_x = q^{-2} \nu_x y,
\\
&& z \nu_x = q^2 \nu_x z, \qquad \qquad z \nu_y = q^{-2} \nu_y z.
\end{eqnarray*}
\end{lemma}

\noindent We are going to show that
$U_{even}$ (resp. 
$U'_{even}$)
is generated by
 $\nu_x,\nu_y,\nu_z,y^{-2}$ (resp.
 $\nu_x,\nu_y,\nu_z$).

\begin{lemma}
\label{lem:six} The following relations hold in
$U_q({\mathfrak{sl}_2})$:
\begin{eqnarray}
&& xy = 1-q^{-1}\nu_z, \qquad \qquad yx = 1-q\nu_z,
\label{eq:xy}
\\
&& yz = 1-q^{-1}\nu_x, \qquad \qquad zy = 1-q\nu_x,
\label{eq:yz}
\\
&& zx = 1-q^{-1}\nu_y, \qquad \qquad xz = 1-q\nu_y.
\label{eq:zx}
\end{eqnarray}
\end{lemma}
\noindent {\it Proof:}  These relations are reformulations of
(\ref{eq:nux})--(\ref{eq:nuz}).
\hfill $\Box$ \\

\begin{lemma}
\label{lem:x2y2z2}
{\rm \cite[Lemma~3.10]{uawe}}.
The following relations hold in 
$U_q({\mathfrak{sl}_2})$:
\begin{eqnarray*}
&&x^2 = 1 - \frac{q\nu_y \nu_z -q^{-1}\nu_z \nu_y}{q-q^{-1}},
\\
&&y^2 = 1 - \frac{q\nu_z \nu_x -q^{-1}\nu_x \nu_z}{q-q^{-1}},
\\
&&z^2 = 1 - \frac{q\nu_x \nu_y -q^{-1}\nu_y \nu_x}{q-q^{-1}}.
\end{eqnarray*}
\end{lemma}

\begin{proposition}
\label{prop:gensets}
The following {\rm (i), (ii)} hold.
\begin{enumerate}
\item[\rm (i)]
The $\F$-algebra
$U_{even}$ is generated by $\nu_x, \nu_y, \nu_z, y^{-2}$;
\item[\rm (ii)]
The $\F$-algebra
$U'_{even}$ is generated by $\nu_x, \nu_y, \nu_z$.
\end{enumerate}
\end{proposition}
\noindent {\it Proof:} 
(i)  Let $W$ denote the subalgebra of
$U_q({\mathfrak{sl}_2})$ generated by
$\nu_x,\nu_y,\nu_z,y^{-2}$.
We show that $W=U_{even}$. Certainly
$W\subseteq U_{even}$, since
each of $\nu_x,\nu_y,\nu_z,y^{-2}$ is contained in
$U_{even}$ by Definition
\ref{def:evenodd} and
(\ref{eq:nux})--(\ref{eq:nuz}).
We now show
$W\supseteq U_{even}$.
By Definition
\ref{def:evenodd}, $U_{even}$ is spanned by
the even elements in the basis
(\ref{eq:xyzbasis}).
Let $x^ry^sz^t$ denote an even element in this
basis. Then $r+s+t=2n$ is even.
Write 
$x^ry^sz^t=g_1g_2\cdots g_n$ 
such that $g_i$ is among
\begin{eqnarray}
\label{eq:list}
x^2, \quad xy, \quad xy^{-1}, \quad xz, \quad y^2, \quad
y^{-2}, \quad yz, \quad y^{-1}z,\quad z^2
\end{eqnarray}
for $1 \leq i \leq n$.
Note that 
$xy^{-1} = xy y^{-2}$
and $y^{-1}z =  y^{-2}yz$.
Now by Lemmas
\ref{lem:six},
\ref{lem:x2y2z2} we see that
$W$ contains each term in (\ref{eq:list}).
Therefore
$W$ contains $g_i$
for $1 \leq i \leq n$, so
 $W$ contains $x^ry^sz^t$. Consequently
$W\supseteq U_{even}$, so
$W=U_{even}$.
\\
\noindent (ii) Similar to the proof of (i).
\hfill $\Box$ \\

\noindent We now consider how the automorphisms $\sigma_a$ act on
$\nu_x$, $\nu_y$, $\nu_z$. 
\begin{lemma}
\label{lem:sigmanu}
For all nonzero $a \in \F$ the automorphism 
$\sigma_a$ sends
\begin{eqnarray*}
\nu_x \mapsto a^2 \nu_x, \quad \qquad
\nu_y \mapsto \nu_y + (a-a^{-1})(aq^2\nu_x - a^{-1}q^{-2}\nu_z)y^{-2},
\qquad \quad 
\nu_z \mapsto a^{-2} \nu_z.
\end{eqnarray*}
\end{lemma}
\noindent {\it Proof:} 
To get the action of $\sigma_a$ on
$\nu_x$ and $\nu_z$, use
Lemma
\ref{lem:sigmasend}
and
(\ref{eq:xzef}).
Now consider the action of $\sigma_a$ on
$\nu_y$. By Lemma
\ref{lem:casequit}
and Lemma
\ref{lem:six},
\begin{eqnarray*}
\Lambda y &=& qzy + q^{-1}xy + q(1-zx)y^2
\\
&=& q(1-q\nu_x) + q^{-1}(1-q^{-1}\nu_z) + \nu_y y^2.
\end{eqnarray*}
In this equation apply $\sigma_a$ to each side and 
use the fact that $\sigma_a$ fixes each of $\Lambda$, $y$.
\hfill $\Box$ \\

\noindent We now consider how the antiautomorphism $\tau$
acts on 
$\nu_x$, $\nu_y$, $\nu_z$. 

\begin{lemma} 
\label{lem:taunu}
The antiautomorphism $\tau$
sends
\begin{eqnarray*}
\nu_x\mapsto \nu_z, \qquad \qquad  
\nu_y\mapsto \nu_y, \qquad \qquad  \nu_z\mapsto \nu_x. 
\end{eqnarray*}
\end{lemma}
\noindent {\it Proof:} 
Use
Lemma
\ref{lem:tau} and
(\ref{eq:nux})--(\ref{eq:nuz}).
\hfill $\Box$ \\

\section{A variation on the Chevalley presentation
of 
$U_q({\mathfrak{sl}_2})$}

\noindent 
In the previous section we saw that 
$U_q({\mathfrak{sl}_2})$ is not generated by
$\nu_x, \nu_y,\nu_z$.
In this section we show that
$U_q({\mathfrak{sl}_2})$ is generated by
$\nu_x, y^{\pm 1},\nu_z$.
Using these generators, we give a presentation of
$U_q({\mathfrak{sl}_2})$ by generators and relations.
As we will see, this presentation is quite close
to the Chevalley presentation.

\begin{lemma} 
\label{lem:nugen}
The elements $\nu_x, y^{\pm 1}, \nu_z$ together
generate 
$U_q({\mathfrak{sl}_2})$. Moreover
\begin{eqnarray}
&&x = y^{-1} - qy^{-1} \nu_z = 
 y^{-1} - q^{-1}\nu_z y^{-1},
\label{eq:g1}
\\
&&z = y^{-1} - q\nu_x y^{-1} =
 y^{-1} - q^{-1}y^{-1}\nu_x.
\label{eq:g2}
\end{eqnarray}
\end{lemma}
\noindent {\it Proof:}  The relations 
(\ref{eq:g1}),
(\ref{eq:g2}) follow from
(\ref{eq:xy}), (\ref{eq:yz}).
The first assertion follows from
(\ref{eq:g1}),
(\ref{eq:g2}) and since
$x,y^{\pm 1},z$ generate
$U_q({\mathfrak{sl}_2})$. 
\hfill $\Box$ \\

\begin{theorem}
\label{thm:nxnzy}
The $\F$-algebra 
$U_q({\mathfrak{sl}_2})$ has a presentation by
generators $\nu_x$, $y^{\pm 1}$, $\nu_z$ and relations
\begin{eqnarray}
&&
\qquad \qquad \quad y y^{-1} = y^{-1} y = 1,
\nonumber
\\
&&
y \nu_x y^{-1} = q^{-2} \nu_x, \qquad \qquad
y \nu_z y^{-1} = q^2 \nu_z,
\label{eq:ynx}
\\
&& \qquad \quad
\frac{q\nu_z \nu_x - q^{-1} \nu_x \nu_z}{q-q^{-1}} = 1-y^2.
\label{eq:nznx}
\end{eqnarray}
An isomorphism with the presentation in
Definition
\ref{def:uqchev} 
sends
\begin{eqnarray*}
\nu_x &\mapsto& -q(q-q^{-1})kf,
\\
y^{\pm 1} &\mapsto & k^{\pm 1},
\\
\nu_z &\mapsto & (q-q^{-1})e.
\end{eqnarray*}
The inverse isomorphism sends
\begin{eqnarray*}
e &\mapsto & (q-q^{-1})^{-1}\nu_z,
\\
k^{\pm 1} &\mapsto & y^{\pm 1},
\\
f &\mapsto & -q^{-1}(q-q^{-1})^{-1}y^{-1}\nu_x.
\end{eqnarray*}

\end{theorem}
\noindent {\it Proof:} 
One checks that each map above is an $\F$-algebra homomorphism,
and that these maps are inverses. 
Therefore each map is an $\F$-algebra isomorphism.
\hfill $\Box$ \\

\begin{definition}\rm 
The presentation of
$U_q({\mathfrak{sl}_2})$ given in
Theorem
\ref{thm:nxnzy} will be called the
{\it modified Chevalley presentation}.
\end{definition}

\section{The element $A$
of
$U_q({\mathfrak{sl}_2})$}

\noindent In Section 9 we will give
another presentation of
$U_q({\mathfrak{sl}_2})$ by generators and relations.
To prepare for this, we introduce the element $A$
of 
$U_q({\mathfrak{sl}_2})$.
\medskip

\noindent Up until the end of Section
11,
 fix a nonzero
$a \in \F$ such that $a^2 \not=1$. Referring to
Lemma
\ref{lem:sigmasend}
we abbreviate
$\sigma=\sigma_a$.

\begin{definition}
\label{def:XZ}
\rm Let $X$ (resp. $Z$) denote the image of
$x$ (resp. $z$)
under $\sigma $.
 By Lemma
\ref{lem:sigform},
\begin{eqnarray}
X = a^{-2}x + (1-a^{-2})y^{-1},
\qquad \qquad 
Z = a^2z + (1-a^2)y^{-1}.
\label{eq:XZ}
\end{eqnarray}
\end{definition}

\begin{lemma} We have
\begin{eqnarray*}
\frac{qXy-q^{-1}yX}{q-q^{-1}} = 1,
\qquad 
\frac{qyZ-q^{-1}Zy}{q-q^{-1}} = 1,
\qquad 
\frac{qZX-q^{-1}XZ}{q-q^{-1}} = 1.
\end{eqnarray*}
\end{lemma}
\noindent {\it Proof:} 
For each relation in
(\ref{eq:xyz}),
apply 
$\sigma$ to each side.
\hfill $\Box$ \\

\begin{lemma}
\label{lem:yi}
We have
\begin{eqnarray*}
y^{-1} = \frac{aX-a^{-1}x}{a-a^{-1}}
= 
\frac{az-a^{-1}Z}{a-a^{-1}}.
\end{eqnarray*}
\end{lemma}
\noindent {\it Proof:}  Use
(\ref{eq:XZ}).
\hfill $\Box$ \\

\begin{corollary} We have
\begin{eqnarray}
a^{-1}x + az
=
aX+a^{-1}Z. 
\label{eq:xzXZ}
\end{eqnarray}
\end{corollary}
\noindent {\it Proof:}  
Use Lemma
\ref{lem:yi}.
\hfill $\Box$ \\

\begin{definition}
\label{def:A} 
\rm
Define
\begin{eqnarray}
A = 
a^{-1}x + az
=
aX+a^{-1}Z. 
\label{eq:A}
\end{eqnarray}
\end{definition}

\noindent Our next general goal is to show that $A,y^{\pm 1}$ together generate
$U_q({\mathfrak{sl}_2})$.

\begin{lemma} 
\label{lem:Ay}
We have
\begin{eqnarray}
Ay &=&
a(1-q\nu_x)
+ 
a^{-1}(1-q^{-1}\nu_z),
\label{eq:Ay}
\\
yA &=&
a(1-q^{-1}\nu_x) +
a^{-1}(1-q\nu_z).
\label{eq:yA}
\end{eqnarray}
\end{lemma}
\noindent {\it Proof:} Use $A=a^{-1}x+az$ and
Lemma
\ref{lem:six}.
\hfill $\Box$ \\

\begin{lemma}\label{lem:nxnz}
We have
\begin{eqnarray}
\nu_x &=& a^{-1}\,\frac{a+a^{-1}}{q+q^{-1}} - 
a^{-1} \frac{qAy-q^{-1}yA}{q^2-q^{-2}},
\label{eq:nxAy}
\\
\nu_z &=& a\,\frac{a+a^{-1}}{q+q^{-1}} - 
a \frac{qyA-q^{-1}Ay}{q^2-q^{-2}}.
\label{eq:nzAy}
\end{eqnarray}
\end{lemma}
\noindent {\it Proof:}  
Solve the equations 
(\ref{eq:Ay}),
 (\ref{eq:yA})
for $\nu_x$, $\nu_z$.
\hfill $\Box$ \\

\begin{lemma} 
\label{lem:Aygen}
The elements $A,y^{\pm 1}$ together generate
$U_q({\mathfrak{sl}_2})$.
\end{lemma}
\noindent {\it Proof:} 
By Lemmas 
\ref{lem:nugen}
and
\ref{lem:nxnz}.
 \hfill $\Box$ \\

\noindent We mention a variation on
Lemma \ref{lem:Ay}.

\begin{lemma} 
\label{lem:Ayconj}
We have
\begin{eqnarray*}
yAy^{-1} &=& a^{-1}q^2 x + a q^{-2}z + (q-q^{-1})(aq^{-1}-a^{-1}q)y^{-1},
\\
y^{-1}Ay &=& a^{-1}q^{-2} x + a q^{2}z + (q-q^{-1})(a^{-1}q^{-1}-aq)y^{-1}.
\end{eqnarray*}
\end{lemma}
\noindent {\it Proof:} 
To verify each equation, eliminate $A$ using
$A=a^{-1}x+az$, and simplify the result using
Lemma
\ref{lem:yconj}.
 \hfill $\Box$ \\

\begin{lemma}  We have
\begin{eqnarray}
yAy^{-1} - (q^2+q^{-2})A + y^{-1}Ay =
-(a+a^{-1})(q-q^{-1})^2 y^{-1}.
\label{lem:yay}
\end{eqnarray}
\end{lemma}
\noindent {\it Proof:} 
To verify 
(\ref{lem:yay}),
evaluate
the left-hand side
using  Lemma
\ref{lem:Ayconj}.
 \hfill $\Box$ \\

\begin{corollary}
\label{cor:yAy}
We have
\begin{eqnarray}
\label{eq:yyA}
y^2 A - (q^2+q^{-2})yAy + A y^2 = -(a+a^{-1})(q-q^{-1})^2y.
\end{eqnarray}
\end{corollary}
\noindent {\it Proof:} 
In equation (\ref{lem:yay})
multiply each term on the left and right by $y$.
 \hfill $\Box$ \\

\noindent We now give $x,z$ in terms of $A,y^{\pm 1}$.

\begin{lemma}
\label{lem:genexpand}
We have
\begin{eqnarray*}
x &=&
A \frac{aq^{-2}}{q^{-2}-q^2} 
+
yAy^{-1} \frac{a}{q^2-q^{-2}}
+
y^{-1} \frac{a(a^{-1}q-aq^{-1})}{q+q^{-1}}
\\
&=&
A \frac{aq^{2}}{q^{2}-q^{-2}} 
+
y^{-1}Ay \frac{a}{q^{-2}-q^{2}}
+
y^{-1} \frac{a(a^{-1}q^{-1}-aq)}{q+q^{-1}}
\end{eqnarray*}
and
\begin{eqnarray*}
z &=& A \frac{a^{-1}q^{2}}{q^{2}-q^{-2}} 
+
yAy^{-1} \frac{a^{-1}}{q^{-2}-q^{2}}
+
y^{-1} \frac{a^{-1}(aq^{-1}-a^{-1}q)}{q+q^{-1}}
\\
 &=& A \frac{a^{-1}q^{-2}}{q^{-2}-q^{2}} 
+
y^{-1}Ay \frac{a^{-1}}{q^{2}-q^{-2}}
+
y^{-1} \frac{a^{-1}(aq-a^{-1}q^{-1})}{q+q^{-1}}.
\end{eqnarray*}
\end{lemma}
\noindent {\it Proof:} 
To verify these equations,
eliminate $yAy^{-1}$ and $y^{-1}Ay$  
using 
Lemma
\ref{lem:Ayconj}. 
 \hfill $\Box$ \\

\noindent
We now describe  how the automorphism $\sigma$ acts on
  $A,y^{\pm 1}$.
\begin{lemma} 
The automorphism $\sigma $ sends
\begin{eqnarray*}
A &\mapsto &  A\frac{a^{-2}q^2-a^2q^{-2}}{q^2-q^{-2}}
+ 
y^{-1}Ay \frac{a^2-a^{-2}}{q^2-q^{-2}}
+
y^{-1}\frac{(a^2-a^{-2})(a^{-1}q-aq^{-1})}{q+q^{-1}},
\\
y^{\pm 1} &\mapsto & y^{\pm 1}.
\end{eqnarray*}
\end{lemma}
\noindent {\it Proof:}  
By construction $\sigma$ fixes $y^{\pm 1}$.
To obtain the image of $A$ under $\sigma$,
in the equation $A=a^{-1}x+az$  first eliminate
$x,z$ using
(\ref{eq:g1}), 
(\ref{eq:g2}). In the resulting equation
apply $\sigma$ to each term and evaluate
using Lemmas
\ref{lem:sigmanu},
\ref{lem:nxnz}.
\hfill $\Box$ \\

\noindent
We now describe  how the antiautomorphism $\tau$ acts on
 $A,y^{\pm 1}$.
\begin{lemma} 
The antiautomorphism $\tau $ sends
\begin{eqnarray*}
A &\mapsto &  A\frac{a^{-2}q^2-a^2q^{-2}}{q^2-q^{-2}}
+ 
yAy^{-1} \frac{a^2-a^{-2}}{q^2-q^{-2}}
+
y^{-1}\frac{(a^2-a^{-2})(a^{-1}q-aq^{-1})}{q+q^{-1}},
\\
y^{\pm 1} &\mapsto & y^{\pm 1}.
\end{eqnarray*}
\end{lemma}
\noindent {\it Proof:}  
Consider the image of $A$ under $\tau$.
Recall that $A=a^{-1}x+az$, and that $\tau$ swaps $x,z$.
Therefore $\tau$ sends
$A\mapsto ax+a^{-1}z$. Evaluating $ax+a^{-1}z$ using
Lemma
\ref{lem:genexpand} we find that the image of $A$ under $\tau$ is
as claimed. Our assertion about $y$ is clear.
\hfill $\Box$ \\

\begin{lemma} 
\label{lem:CasAy}
The Casimir element $\Lambda$ looks as follows
in terms of $A, y^{\pm 1}$:
\begin{eqnarray*}
\Lambda = \frac{
A^2y-(q^2+q^{-2})AyA+yA^2+(q^2-q^{-2})^2y+(a+a^{-1})(q-q^{-1})^2 A}
{(q-q^{-1})(q^2-q^{-2})}.
\end{eqnarray*}
\end{lemma}
\noindent {\it Proof:} 
Using $A=a^{-1}x+az$ we find that $A^2y-(q^2+q^{-2})AyA+yA^2$
is equal to $a^{-2}$ times
\begin{eqnarray}
x^2y-(q^2+q^{-2})xyx+yx^2 
\label{eq:xgxy}
\end{eqnarray}
plus $a^2$ times
\begin{eqnarray}
z^2y-(q^2+q^{-2})zyz+yz^2 
\label{eq:zzy}
\end{eqnarray}
plus
\begin{eqnarray}
yzx + zxy - (q^2+q^{-2})xyz 
+ xzy + yxz - (q^2+q^{-2})zyx.
\label{eq:six}
\end{eqnarray}
By Lemma
\ref{lem:qsymrel}, the expressions 
(\ref{eq:xgxy}) and
(\ref{eq:zzy}) are equal to
 $-(q-q^{-1})^2 x$ and
 $-(q-q^{-1})^2 z$, respectively.
By Lemma
\ref{lem:casequit} the expression
(\ref{eq:six}) is equal to
\begin{eqnarray*}
\Lambda (q-q^{-1})(q^2-q^{-2})
- x(q-q^{-1})^2 
- y(q^2-q^{-2})^2 
- z(q-q^{-1})^2. 
\end{eqnarray*}
The result follows from these comments after a brief calculation.
\hfill $\Box$ \\

\noindent We mention several variations on Lemma
\ref{lem:CasAy}.

\begin{lemma}\label{lem:Lamfromnxnz}
We have
\begin{eqnarray}
&&\Lambda = a\frac{q\nu_x A - q^{-1} A \nu_x}{q-q^{-1}} + (q+q^{-1})y,
\label{eq:Lamfromnx}
\\
&&\Lambda = a^{-1}\frac{q A\nu_z - q^{-1} \nu_z A}{q-q^{-1}} + (q+q^{-1})y.
\label{eq:Lamfromnz}
\end{eqnarray}
\end{lemma}
\noindent {\it Proof:}
To verify 
(\ref{eq:Lamfromnx}), eliminate $\nu_x$ using
(\ref{eq:nxAy}) and compare the result with Lemma
\ref{lem:CasAy}.
Equation 
(\ref{eq:Lamfromnz}) is similarly verified, by eliminating $\nu_z$ using
(\ref{eq:nzAy}).
\hfill $\Box$ \\

\noindent In Corollary
\ref{cor:yAy} we gave an equation relating $A,y$. We now
give some more equations relating $A,y$.
 Recall the notation
\begin{eqnarray*}
\lbrack n \rbrack_q = \frac{q^n-q^{-n}}{q-q^{-1}}
\qquad \qquad n \in \Z.
\end{eqnarray*}

\begin{proposition} 
\label{prop:tdrel}
The elements $A$, $y$ satisfy
\begin{eqnarray}
&&
y^3A-\lbrack 3 \rbrack_q y^2 Ay+ \lbrack 3 \rbrack_q yAy^2-Ay^3 = 0,
\label{eq:td1}
\\
&&
\frac{A^3y-\lbrack 3 \rbrack_q A^2yA+\lbrack 3 \rbrack_q AyA^2-yA^3}
{(q^2-q^{-2})^2} = 
 yA-Ay.
\label{eq:td2}
\end{eqnarray}
\end{proposition}
\noindent {\it Proof:}
To get
(\ref{eq:td1}),
take the commutator of $y$ with each side of
(\ref{eq:yyA}).
To get 
(\ref{eq:td2}),
take the commutator
of $A$ with each side of
the equation in  
 Lemma
\ref{lem:CasAy}.
\hfill $\Box$ \\

\begin{note}\rm The equations 
(\ref{eq:td1}), 
(\ref{eq:td2}) 
are a special case of
the tridiagonal relations
\cite{td99}.
\end{note}

\begin{note}\rm In
\cite[Section~23]{boyd}
C.~Worawannotai shows how
the elements $A,y$ are related to Leonard pairs of
dual $q$-Krawtchouk type. 
\end{note}


\begin{lemma}
\label{lem:yAy}
The elements $A,y$ satisfy
\begin{eqnarray*}
&&Ay^2A+AyAy+yAyA-\lbrack 3 \rbrack_q yA^2y 
+(a+a^{-1})(q-q^{-1})^2(Ay+yA)
\\
&& \qquad = \;
(q^2-q^{-2})^2 y^2
+ (aq-a^{-1}q^{-1})(aq^{-1}-a^{-1}q)(q-q^{-1})^2.
\end{eqnarray*}
\end{lemma}
\noindent {\it Proof:}  
In the equation (\ref{eq:nznx}) 
eliminate $\nu_x$, $\nu_z$ using
 Lemma
\ref{lem:nxnz} and simplify the result.
\hfill $\Box$ \\

\section{The antiautomorphism $\dagger$}

\noindent In this section we introduce the
antiautomorphism $\dagger$ of 
$U_q({\mathfrak{sl}_2})$, and describe it 
from several points of view.

\begin{definition}
\label{def:dagger}
\rm Define a map
$
\dagger: 
U_q({\mathfrak{sl}_2}) \to
U_q({\mathfrak{sl}_2})$
to be the composition
\[
\dagger: \quad 
\begin{CD}
 U_{q}(\mathfrak{sl}_2)  @>>\tau >
  U_q(\mathfrak{sl}_2)  @>>\sigma> 
  U_q(\mathfrak{sl}_2) 
		                     \end{CD} 
				                   \]
where $\tau$
is from
Lemma \ref{lem:tau} and
$\sigma$ is 
from above Definition
\ref{def:XZ}.

\end{definition}

\begin{lemma} 
\label{lem:xZzX}
The map $\dagger$ is the unique antiautomorphism
of
  $U_q(\mathfrak{sl}_2) $ that sends
  \begin{eqnarray}
x \mapsto Z, \qquad \qquad y^{\pm 1} \mapsto y^{\pm 1}, 
\qquad \qquad z \mapsto X.
\label{eq:dagger}
\end{eqnarray}
\end{lemma}
\noindent {\it Proof:}  
The map $\dagger$ is an antiautomorphism since
$\tau$ is an antiautomorphism and $\sigma$ is an
automorphism.
The antiautomorphism $\dagger$ satisfies (\ref{eq:dagger}) 
by Lemma
\ref{lem:tau} and
Definition
\ref{def:XZ}.
The antiautomorphism $\dagger$ is unique since
$x,y^{\pm 1},z$ generate 
  $U_q(\mathfrak{sl}_2) $.
\hfill $\Box$ \\

\begin{lemma}
\label{lem:nxnzdag}
The map $\dagger$ is the
unique antiautomorphism of 
  $U_q(\mathfrak{sl}_2) $ that sends
\begin{eqnarray}
\label{eq:nxdagsend}
\nu_x \mapsto a^{-2} \nu_z, \qquad \qquad
y^{\pm 1}\mapsto y^{\pm 1},
\qquad \qquad
\nu_z \mapsto a^{2} \nu_x.
\end{eqnarray}
\end{lemma}
\noindent {\it Proof:} By
Lemma
\ref{lem:xZzX} the
 map
$\dagger$ is an antiautomorphism and
fixes $y^{\pm 1}$.
The map $\dagger$ satisfies the rest of 
(\ref{eq:nxdagsend}) by
Lemmas 
\ref{lem:sigmanu},
\ref{lem:taunu}
and
Definition
\ref{def:dagger}.
The antiautomorphism $\dagger$ is unique since
$\nu_x, y^{\pm 1}, \nu_z$ generate 
  $U_q(\mathfrak{sl}_2) $ by
Lemma
\ref{lem:nugen}.
\hfill $\Box$ \\

\begin{lemma} 
\label{lem:fixAy}
The map $\dagger$ is the
unique antiautomorphism of 
  $U_q(\mathfrak{sl}_2) $ that fixes each of
  $A,y^{\pm 1}$.
\end{lemma}
\noindent {\it Proof:} By
Lemma
\ref{lem:xZzX} the
 map
$\dagger$ is an antiautomorphism and
fixes $y^{\pm 1}$.
By
Definition
\ref{def:A} and
Lemma
\ref{lem:xZzX},
the image of
$A=a^{-1}x+az$ under $\dagger$ is
$a^{-1}Z+aX=A$.
The antiautomorphism $\dagger$ is unique since
$A, y^{\pm 1}$ generate 
  $U_q(\mathfrak{sl}_2) $ by
Lemma
\ref{lem:Aygen}.
\hfill $\Box$ \\

\begin{lemma} The antiautomorphism $\dagger$ fixes $\Lambda$.
\end{lemma}
\noindent {\it Proof:}  By Definition
\ref{def:dagger} and since each of $\tau, \sigma$ fixes $\Lambda$.
\hfill $\Box$ \\

\begin{lemma}
\label{lem:daginv}
We have $\dagger^2=1$.
\end{lemma}
\noindent {\it Proof:} 
Use Lemma
\ref{lem:tauprod}
and
Definition
\ref{def:dagger}.
\hfill $\Box$ \\

\begin{lemma} 
\label{lem:factor}
We have the factorizations
\[
\sigma: \quad 
\begin{CD}
 U_{q}(\mathfrak{sl}_2)  @>>\tau >
  U_q(\mathfrak{sl}_2)  @>>\dagger> 
  U_q(\mathfrak{sl}_2) 
		                     \end{CD},
				                   \]
\[
\sigma^{-1}: \quad 
\begin{CD}
 U_{q}(\mathfrak{sl}_2)  @>>\dagger >
  U_q(\mathfrak{sl}_2)  @>>\tau> 
  U_q(\mathfrak{sl}_2) 
		                     \end{CD}.
				                   \]
\end{lemma}
\noindent {\it Proof:} 
Use
Definition
\ref{def:dagger} along with
Lemma \ref{lem:daginv} and the last assertion
of Lemma
\ref{lem:tau}.
\hfill $\Box$ \\

\noindent In Lemma
\ref{lem:genexpand} we found each of $x,z$ in terms
of $A,y^{\pm 1}$. We now do something similar for
$X,Z$.

\begin{lemma} The elements $X$, $Z$ look as follows in terms
of $A, y^{\pm 1}$:
\begin{eqnarray*}
X &=&
A \frac{a^{-1}q^{-2}}{q^{-2}-q^2} 
+
yAy^{-1} \frac{a^{-1}}{q^2-q^{-2}}
+
y^{-1} \frac{a^{-1}(aq-a^{-1}q^{-1})}{q+q^{-1}}
\\
&=&
A \frac{a^{-1}q^{2}}{q^{2}-q^{-2}} 
+
y^{-1}Ay \frac{a^{-1}}{q^{-2}-q^{2}}
+
y^{-1} \frac{a^{-1}(aq^{-1}-a^{-1}q)}{q+q^{-1}}
\end{eqnarray*}
and
\begin{eqnarray*}
Z &=& A \frac{aq^{2}}{q^{2}-q^{-2}} 
+
yAy^{-1} \frac{a}{q^{-2}-q^{2}}
+
y^{-1} \frac{a(a^{-1}q^{-1}-aq)}{q+q^{-1}}
\\
 &=& A \frac{aq^{-2}}{q^{-2}-q^{2}} 
+
y^{-1}Ay \frac{a}{q^{2}-q^{-2}}
+
y^{-1} \frac{a(a^{-1}q-aq^{-1})}{q+q^{-1}}.
\end{eqnarray*}
\end{lemma}
\noindent {\it Proof:}  
In each of the four equations of
Lemma
\ref{lem:genexpand},
apply $\dagger$ 
to each side
and evaluate the results using
Lemmas
\ref{lem:xZzX},
\ref{lem:fixAy}.
\hfill $\Box$ \\

\section{The invariant presentation of 
$U_q({\mathfrak{sl}_2})$}

\noindent In this section we give a presentation
of
$U_q({\mathfrak{sl}_2})$ by generators and relations,
using the generators $A$, $y^{\pm 1}$.

\begin{theorem}
\label{thm:invar}
The $\F$-algebra 
$U_q({\mathfrak{sl}_2})$ has a presentation by generators
$A,y^{\pm 1}$ and relations
\begin{eqnarray*}
&& 
\qquad \qquad \qquad \qquad \qquad
yy^{-1} = y^{-1}y = 1,
\\
&&\qquad y^2A-(q^2+q^{-2})yAy+Ay^2= -(a+a^{-1})(q-q^{-1})^2y,
\\
&&Ay^2A+AyAy+yAyA-\lbrack 3 \rbrack_q yA^2y + 
(a+a^{-1})(q-q^{-1})^2(Ay+yA)
\\
&& \qquad = \; 
(q^2-q^{-2})^2 y^2 +
(aq-a^{-1}q^{-1})(aq^{-1}-a^{-1}q)(q-q^{-1})^2.
\end{eqnarray*}
\noindent An isomorphism with the presentation in
Theorem \ref{thm:nxnzy}
 sends
\begin{eqnarray}
A &\mapsto & a(1-q\nu_x)y^{-1}+ a^{-1}(1-q^{-1}\nu_z)y^{-1},
\label{eq:Asend}
\\
y^{\pm 1} &\mapsto & y^{\pm 1}.
\label{eq:ysend}
\end{eqnarray}
\noindent The inverse isomorphism sends
\begin{eqnarray}
\nu_x &\mapsto& a^{-1}\,\frac{a+a^{-1}}{q+q^{-1}} - 
a^{-1} \frac{qAy-q^{-1}yA}{q^2-q^{-2}},
\label{eq:nxsend}
\\
y^{\pm 1} &\mapsto & y^{\pm 1},
\label{eq1:ysend}
\\
\nu_z &\mapsto& a\,\frac{a+a^{-1}}{q+q^{-1}} - 
a \frac{qyA-q^{-1}Ay}{q^2-q^{-2}}.
\label{eq:nzsend}
\end{eqnarray}
\end{theorem}
\noindent {\it Proof:}
The map
(\ref{eq:Asend}),
(\ref{eq:ysend})
gives an $\F$-algebra homomorphism 
by Corollary 
\ref{cor:yAy} and Lemma
\ref{lem:yAy}.
One checks that the map
(\ref{eq:nxsend})--(\ref{eq:nzsend})
gives an $\F$-algebra homomorphism.
One also checks that 
the above maps are inverses.
Therefore each map is an $\F$-algebra isomorphism.
\hfill $\Box$ \\

\begin{note}
\label{note:invar}
\rm 
For the presentation of
$U_q({\mathfrak{sl}_2})$ given in Theorem
\ref{thm:invar}, the relations
are invariant 
under each of the following moves:
(i) $q\mapsto q^{-1}$;
(ii) $a\mapsto a^{-1}$.
\end{note}

\begin{definition}\rm In view of
Lemma
\ref{lem:fixAy}
and
Note \ref{note:invar},
the presentation for
$U_q({\mathfrak{sl}_2})$ given in Theorem
\ref{thm:invar} will be called the 
{\it invariant presentation}.
\end{definition}

\section{The algebra
$U^{\vee}_q$}

\noindent We turn our attention to a certain
subalgebra of 
 $U_q({\mathfrak{sl}_2})$, which is denoted
$U^{\vee}_q$. We will define
$U^{\vee}_q$ after some motivational comments.
Recall the
 equitable generators
$x,y^{\pm 1},z$ of 
$U_q({\mathfrak{sl}_2})$
from Section 3,
 and the elements $X,Z,A$ from Section 7.
Consider the subspace of
$U_q({\mathfrak{sl}_2})$
spanned by
$x,y^{-1},z$. This subspace contains 
$X,Z, A$  by 
(\ref{eq:XZ}),
(\ref{eq:A}), and is invariant under
the maps $\sigma, \tau, \dagger$
by
Lemmas
\ref{lem:sigform},
\ref{lem:tau}
and Definition
\ref{def:dagger}. 

\begin{definition}
\label{def:Uvee}
\rm
Let 
$U^{\vee}_q$ denote the subalgebra of
$U_q({\mathfrak{sl}_2})$ generated by
$x$, $y^{-1}$, $z$.
\end{definition}
\noindent Our next goal is to
(i) display a basis for $U^\vee_q$;
(ii) give several presentations of $U^\vee_q$ 
by generators and relations.

\begin{lemma}  The following relations hold in
$U^{\vee}_q$:
\begin{eqnarray}
\frac{qy^{-1}x-q^{-1}xy^{-1}}{q-q^{-1}} = y^{-2},
\qquad 
\frac{qzy^{-1}-q^{-1}y^{-1}z}{q-q^{-1}} = y^{-2},
\qquad 
\frac{qzx-q^{-1}xz}{q-q^{-1}} = 1.
\label{eq1:Uveepres}
\end{eqnarray}
\end{lemma}
\noindent {\it Proof:} 
Adjust the relations
(\ref{eq:xyz}).
\hfill $\Box$ \\

\noindent 
As we investigate
$U^{\vee}_q$ it is useful to consider the following
algebra.

\begin{definition}
\label{def:lineuq}
\rm Let $\overline U_q$ denote
the $\F$-algebra defined by generators $x$, $\overline y$, $z$
and relations
\begin{eqnarray}
\frac{q\overline y x-q^{-1}x\overline y}{q-q^{-1}} = \overline y^{2},
\qquad 
\frac{qz\overline y-q^{-1}\overline yz}{q-q^{-1}} = \overline y^{2},
\qquad 
\frac{qzx-q^{-1}xz}{q-q^{-1}} = 1.
\label{eq:hlineUpres}
\end{eqnarray}
\end{definition}

\begin{lemma} \label{lem:barvee}
There exists an $\F$-algebra homomorphism 
$\iota: \overline U_q \to 
U_q({\mathfrak{sl}_2})$
 that sends
\begin{eqnarray*}
x \mapsto x, \qquad \qquad \overline y \mapsto y^{-1}, \qquad \qquad 
z \mapsto z.
\end{eqnarray*}
\end{lemma}
\noindent {\it Proof:}  
Compare 
(\ref{eq1:Uveepres}) and
(\ref{eq:hlineUpres}).
\hfill $\Box$ \\

\noindent Note that $U^\vee_q$ is the image of
$\overline U_q$ under $\iota$.

\begin{lemma}
\label{lem:Uveebasis} 
The following is a basis for the
$\F$-vector space $\overline U_q$:
\begin{eqnarray}
x^r \overline y^s z^t \qquad \qquad r,s,t \in \N.
\label{eq:Uveespan}
\end{eqnarray}
\end{lemma}
\noindent {\it Proof:} 
The vectors 
(\ref{eq:Uveespan})
are linearly independent, since their
images under $\iota $ are linearly
independent by
Lemma
\ref{lem:Uqbasis}.
We show that the vectors
(\ref{eq:Uveespan}) span
$\overline U_q$. 
For $n \in \N$, by a {\it word of length $n$}
in $\overline U_q$
we mean a product
$x_1x_2 \cdots x_n$ such that
$x_i$ is among $x,\overline y,z$ for
$1 \leq i \leq n$. We interpret the word of length 0
to be the multiplicative identity of
$\overline U_q$. By definition 
$\overline U_q$ is spanned by the
words.
Pick a word $w$ and write $w=x_1x_2 \cdots x_n$.
By the {\it $(x,z)$-length} of $w$
we mean the number of elements among
$x_1, x_2, \ldots, x_n$ that are equal to $x$ or $z$.
By an {\it inversion} in $w$
we mean an ordered pair of elements  $(x_i,x_j)$
$(1 \leq i,j\leq n)$ such that
$i<j$  and $(x_i,x_j)$ is 
$(\overline y,x)$ or $(z, \overline y)$ or $(z,x)$.
The word $w$
is called {\it reducible} whenever it
has at least one inversion, and {\it irreducible} otherwise.
The list (\ref{eq:Uveespan}) consists of the irreducible words.
The words (\ref{eq:Uveespan}) span a subspace of
$\overline U_q$ that we denote by $L$.
We show
$L=\overline U_q $.
In order to do this, it suffices
to show that each word is contained in
$L$.  Suppose there exists a word 
that is not contained in $L$.
Let $N$ denote the minimal $(x,z)$-length among all such words.
Consider the set of words that are not 
contained in $L$ and 
have
$(x,z)$-length $N$.
Pick a word $w$ in this set that has a minimal number of
inversions. The word $w$ is reducible;
otherwise $w$ is included in
(\ref{eq:Uveespan}) and hence $w \in L$, for a contradiction.
Write $w=x_1x_2\cdots x_n$.
Since $w$ is reducible, there exists an integer
$i$ $(2 \leq i \leq n)$ such that $(x_{i-1},x_i)$
is an inversion. There are three cases.
First assume that
$(x_{i-1},x_i)=(\overline y,x)$.
Note that
$\overline yx = q^{-2}x\overline y+
(1-q^{-2})\overline y^{2}$.
Therefore $w=q^{-2}u+(1-q^{-2})v$ where
$u =  x_1\cdots x_{i-2}x\overline yx_{i+1}\cdots x_n$ and
$v =  x_1\cdots x_{i-2}\overline y^{2}x_{i+1}\cdots x_n$.
The word $u$ has $(x,z)$-length $N$ and
one fewer inversions than $w$. 
Therefore $u \in L$.
The word $v$ has $(x,z)$-length $N-1$,
so $v \in L$.
Therefore $w \in L$ for 
a contradiction.
Next assume that
$(x_{i-1},x_i)=(z,\overline y)$.
Note that $z\overline y = q^{-2}\overline yz+
(1-q^{-2})\overline y^{2}$.
Therefore $w=q^{-2}u+(1-q^{-2})v$ where
$u =  x_1\cdots x_{i-2}\overline yzx_{i+1}\cdots x_n$ and
$v =  x_1\cdots x_{i-2}\overline y^{2}x_{i+1}\cdots x_n$.
The word $u$ has $(x,z)$-length $N$ and
one fewer inversions than $w$. 
Therefore $u \in L$.
The word $v$ has $(x,z)$-length $N-1$, 
so $v \in L$.
Therefore $w \in L$ for 
a contradiction.
Next assume that
$(x_{i-1},x_i)=(z,x)$.
Note that $zx = q^{-2}xz+1-q^{-2}$.
Therefore $w=q^{-2}u+(1-q^{-2})v$ where
$u =  x_1\cdots x_{i-2}xzx_{i+1}\cdots x_n$ and
$v =  x_1\cdots x_{i-2}x_{i+1}\cdots x_n$.
The word $u$ has $(x,z)$-length $N$ and
one fewer inversions than $w$. 
Therefore $u \in L$.
The word $v$ has $(x,z)$-length $N-2$, so
$v \in L$. Therefore $w \in L$ for 
a contradiction.
In each of the three cases we obtained a contradiction.
Therefore every word is contained in $L$,
and consequently
$L=\overline U_q$. The result follows.
\hfill $\Box$ \\

\begin{proposition}
\label{prop:inj}
The map 
$\iota:\overline U_q \to 
U_q({\mathfrak{sl}_2})$
is injective.
\end{proposition}
\noindent {\it Proof:} 
For the basis vectors
(\ref{eq:Uveespan}) their images under $\iota$
are linearly independent by Lemma
\ref{lem:Uqbasis}.
\hfill $\Box$ \\

\noindent 
Proposition \ref{prop:inj} gives an isomorphism
of $\F$-algebras $\iota :\overline U_q \to U^\vee_q$.

\begin{corollary}
\label{cor:Uqbasisyi}
The following is a basis for the $\F$-vector
space $U^{\vee}_q$:
\begin{eqnarray*}
x^r y^{-s} z^t \qquad \qquad r,s,t \in \N.
\end{eqnarray*}
\end{corollary}
\noindent {\it Proof:}  Apply the isomorphism
$\iota : \overline U_q \to U^\vee_q$ to the
basis
(\ref{eq:Uveespan}).
\hfill $\Box$ \\

\begin{corollary}
\label{cor:uveeuprime}
The following {\rm (i), (ii)} coincide:
\begin{enumerate}
\item[\rm (i)]
the intersection of $U^\vee_q $ and $U'_q$;
\item[\rm (ii)]
the subalgebra of 
$U_q({\mathfrak{sl}_2})$ generated by $x,z$.
\end{enumerate}
This common subalgebra has a basis
\begin{eqnarray*}
x^rz^t \qquad \qquad r,t \in \N.
\end{eqnarray*}
\end{corollary}
\noindent {\it Proof:} Compare the basis for
$U'_q$ given in 
Lemma \ref{lem:uprimebasis}, with the
basis for 
$U^\vee_q$ given in 
Corollary
\ref{cor:Uqbasisyi}.
\hfill $\Box$ \\

\begin{theorem} 
\label{lem:Uveepres}
The $\F$-algebra 
$U^{\vee}_q$ has a presentation by
generators $x,y^{-1},z$ and relations
\begin{eqnarray}
\frac{qy^{-1}x-q^{-1}xy^{-1}}{q-q^{-1}} = y^{-2},
\qquad 
\frac{qzy^{-1}-q^{-1}y^{-1}z}{q-q^{-1}} = y^{-2},
\qquad 
\frac{qzx-q^{-1}xz}{q-q^{-1}} = 1.
\label{eq:Uveepres}
\end{eqnarray}
\end{theorem}
\noindent {\it Proof:} 
By Definition
\ref{def:lineuq} and since
$\iota :\overline U_q \to U^\vee_q$ is
an isomorphism of $\F$-algebras.
\hfill $\Box$ \\

\noindent We now clarify an aspect of
$U^\vee_q$.

\begin{lemma} 
\label{cor:yinotinv} 
The element $y^{-1}$ is not invertible in $U^{\vee}_q$.
\end{lemma}
\noindent {\it Proof:} 
Suppose that
$y^{-1}$ is invertible in
 $U^{\vee}_q$, and let $\xi \in U^\vee_q$ denote its inverse.
 In the algebra 
$U_q({\mathfrak{sl}_2})$,
each of $y,\xi$ is an inverse of $y^{-1}$. The inverse of $y^{-1}$ is
unique, so $y=\xi$.
Comparing 
Lemma
\ref{lem:Uqbasis} and
Corollary
\ref{cor:Uqbasisyi},
we see that
$y \not\in U^{\vee}_q$. This is a contradiction,
and the result follows.
\hfill $\Box$ \\

\begin{corollary} 
\label{cor:byinotinv} 
The element $\overline y$ is not invertible in
$\overline U_q$. 
\end{corollary}
\noindent {\it Proof:} 
By Lemma
\ref{cor:yinotinv} and since
$\iota$ sends
$\overline y\mapsto y^{-1}$.
\hfill $\Box$ \\

\noindent In Theorem
\ref{lem:Uveepres}
we displayed a presentation
of $U^{\vee}_q$ by generators and relations. Shortly we
will display several more.

\begin{definition}\rm Let $S$ denote the subspace of
$U_q({\mathfrak{sl}_2})$ spanned by $x,y^{-1},z$.
\end{definition}

\begin{lemma} 
\label{lem:3gen}
Each of the following {\rm (i)--(iv)} is
a basis for the $\F$-vector space $S$:
\begin{eqnarray*}
{\rm (i)}\;
x,y^{-1},z;
\qquad \qquad
{\rm (ii)}\;
X,y^{-1},Z;
\qquad \qquad
{\rm (iii)}\;
x,X,A;
\qquad \qquad
{\rm (iv)}\;
z,Z,A.
\end{eqnarray*}
\end{lemma}
\noindent {\it Proof:} The vectors (i) are linearly
independent, by Lemma
\ref{lem:Uqbasis} or
Corollary
\ref{cor:Uqbasisyi}.
To finish the proof use
(\ref{eq:XZ})
and $A=a^{-1}x+az$.
\hfill $\Box$ \\

\noindent
We now clarify the relationship between
the four bases for $S$ 
given in Lemma
\ref{lem:3gen}. We will be discussing
transition matrices and matrix representations,
following the conventions of
\cite[Sections~10,~15]{uq}.

\begin{lemma}
\label{lem:3trans}
We refer to the bases {\rm (i)--(iv)} in Lemma
\ref{lem:3gen}.
The transition matrices from {\rm (i)} to {\rm (ii)} and
{\rm (ii)} to {\rm (i)} are, respectively,
\begin{eqnarray*}
\left(
\begin{array}{c c c }
a^{-2} & 0  & 0\\
1-a^{-2} & 1  & 1-a^2 \\
0 & 0 & a^2 
\end{array}
\right),
\qquad \qquad
\left(
\begin{array}{c c c }
a^{2} & 0 & 0\\
1-a^2 & 1  & 1-a^{-2} \\
0 & 0 &  a^{-2} 
\end{array}
\right).
\end{eqnarray*}
The transition matrices from {\rm (i)} to {\rm (iii)} and
{\rm (iii)} to {\rm (i)} are, respectively,
\begin{eqnarray*}
\left(
\begin{array}{c c c }
1 & a^{-2}  & a^{-1}\\
0 & 1-a^{-2}  & 0 \\
0 & 0 & a 
\end{array}
\right),
\qquad \qquad
\left(
\begin{array}{c c c }
1 & \frac{a^{-1}}{a^{-1}-a}  & -a^{-2}\\
0 & \frac{a}{a-a^{-1}}  & 0 \\
0 & 0 & a^{-1} 
\end{array}
\right).
\end{eqnarray*}
The transition matrices from {\rm (i)} to {\rm (iv)} and
{\rm (iv)} to {\rm (i)} are, respectively,
\begin{eqnarray*}
\left(
\begin{array}{c c c }
0 & 0  & a^{-1}\\
0 & 1-a^{2}  & 0 \\
1 & a^2 & a 
\end{array}
\right),
\qquad \qquad
\left(
\begin{array}{c c c }
-a^{2} & \frac{a}{a-a^{-1}}  & 1\\
0 & \frac{a^{-1}}{a^{-1}-a}  & 0 \\
a & 0 &  0 
\end{array}
\right).
\end{eqnarray*}
\end{lemma}
\noindent {\it Proof:}
Use
(\ref{eq:XZ})
and $A=a^{-1}x+az$.
\hfill $\Box$ \\

\noindent By our comments above
Definition
\ref{def:Uvee}, the 
 subspace $S$ is invariant under 
$\tau$ and $\dagger$. We now describe the actions
of $\tau$ and $\dagger$ on $S$.

\begin{lemma} 
\label{lem:Staudag}
We refer to the bases
{\rm (i)--(iv)} in
Lemma
\ref{lem:3gen}.
With respect to  {\rm (i)} the matrices representing
$\tau$ and $\dagger$ are, respectively,
\begin{eqnarray*}
\left(
\begin{array}{c c c }
0 & 0  & 1\\
0 & 1 & 0 \\
1 & 0 & 0 
\end{array}
\right),
\qquad \qquad
\left(
\begin{array}{c c c }
0 & 0  & a^{-2}\\
1-a^2 & 1  & 1-a^{-2} \\
a^2& 0 & 0 
\end{array}
\right).
\end{eqnarray*}
With respect to {\rm (ii)} the matrices representing
$\tau$ and $\dagger$ are, respectively,
\begin{eqnarray*}
\left(
\begin{array}{c c c }
0 & 0  & a^4\\
a^{-2}(a^2-a^{-2}) & 1 & a^2(a^{-2}-a^2) \\
a^{-4} & 0 & 0 
\end{array}
\right),
\qquad \qquad
\left(
\begin{array}{c c c }
0 & 0  & a^2\\
1-a^{-2} & 1 & 1-a^2 \\
a^{-2} & 0 & 0 
\end{array}
\right).
\end{eqnarray*}
With respect to {\rm (iii)} the matrices representing
$\tau$ and $\dagger$ are, respectively,
\begin{eqnarray*}
\left(
\begin{array}{c c c }
-a^{-2} & -a^{-3}(a+a^{-1})  & a^{-1}(a^2-a^{-2})\\
0 & 1 & 0 \\
a^{-1} & a^{-3} & a^{-2} 
\end{array}
\right),
\qquad \qquad
\left(
\begin{array}{c c c }
0 & -a^{-2}  & 0\\
-a^2 & 0 & 0 \\
a & a^{-1} & 1 
\end{array}
\right).
\end{eqnarray*}
With respect to {\rm (iv)} the matrices representing
$\tau$ and $\dagger$ are, respectively,
\begin{eqnarray*}
\left(
\begin{array}{c c c }
-a^{2} & -a^{3}(a+a^{-1})  & a(a^{-2}-a^{2})\\
0 & 1 & 0 \\
a & a^{3} & a^{2} 
\end{array}
\right),
\qquad \qquad
\left(
\begin{array}{c c c }
0 & -a^{2}  & 0\\
-a^{-2} & 0 & 0 \\
a^{-1} & a & 1 
\end{array}
\right).
\end{eqnarray*}
\end{lemma}
\noindent {\it Proof:}
Use Lemmas
\ref{lem:tau},
\ref{lem:xZzX},
\ref{lem:3trans}
and linear algebra.
\hfill $\Box$ \\

\noindent The action of the automorphism
$\sigma$ on $S$ is found using
Lemmas
\ref{lem:factor}
and
\ref{lem:Staudag}. 

\begin{lemma}
The subalgebra
$U^{\vee}_q$ is invariant under each 
map $\tau$, $\dagger$, $\sigma$.
\end{lemma}
\noindent {\it Proof:} The algebra  
$U^{\vee}_q$ is generated by $S$, and $S$ is invariant
under each of
$\tau$, $\dagger$, $\sigma$.
\hfill $\Box$ \\

\noindent 
Consider the sets
{\rm (i)}--{\rm (iv)} 
in Lemma
\ref{lem:3gen}. Each set
 is a generating set
for 
$U^{\vee}_q$. For case (i), the corresponding
presentation of 
$U^{\vee}_q$ by generators and relations was given in
Theorem
\ref{lem:Uveepres}.
For the cases
 (ii)--(iv) we now do something similar.

\begin{theorem}
\label{thm:uvpres2}
The $\F$-algebra $U^{\vee}_q$ has a presentation
by generators $X,y^{-1},Z$ and relations
\begin{eqnarray*}
\frac{qy^{-1}X-q^{-1}Xy^{-1}}{q-q^{-1}} = y^{-2},
\qquad 
\frac{qZy^{-1}-q^{-1}y^{-1}Z}{q-q^{-1}} = y^{-2},
\qquad 
\frac{qZX-q^{-1}XZ}{q-q^{-1}} = 1.
\end{eqnarray*}
\end{theorem}
\noindent {\it Proof:} The restriction of
$\sigma$ to 
$U^{\vee}_q$ is an automorphism of
$U^{\vee}_q$ that sends the generators
$x,y^{-1},z$ to the generators
$X,y^{-1},Z$ respectively.
The result follows from this and
Theorem
\ref{lem:Uveepres}.
\hfill $\Box$ \\

\begin{theorem}
\label{thm:uveepres3}
The $\F$-algebra $U^{\vee}_q$ has a presentation
by generators $x,X,A$ and relations
\begin{eqnarray*}
&&
\frac{qAx-q^{-1}xA}{q-q^{-1}} = a^{-1}x^2+a,
\qquad \qquad
\frac{qAX-q^{-1}XA}{q-q^{-1}} = a X^2+a^{-1},
\\
&&
\qquad
a^{-1}x^2 -\frac{a^{-1}q-aq^{-1}}{q-q^{-1}}xX-
\frac{aq-a^{-1}q^{-1}}{q-q^{-1}}Xx + aX^2 = 0.
\end{eqnarray*}
An isomorphism with the presentation in
Theorem
\ref{lem:Uveepres} sends
\begin{eqnarray*}
x &\mapsto& x,
\\
X &\mapsto& a^{-2}x + (1-a^{-2})y^{-1},
\\
A &\mapsto& a^{-1}x+az.
\end{eqnarray*}
The inverse isomorphism sends
\begin{eqnarray*}
x &\mapsto& x,
\\
y^{-1} &\mapsto& \frac{aX-a^{-1}x}{a-a^{-1}},
\\
z &\mapsto& a^{-1}A-a^{-2}x.
\end{eqnarray*}
\end{theorem}
\noindent {\it Proof:} 
One routinely checks that each map is an $\F$-algebra homomorphism,
and that the two maps are inverses. Therefore each map is
an $\F$-algebra isomorphism.
\hfill $\Box$ \\

\begin{theorem} 
\label{thm:uveepres4}
The $\F$-algebra $U^{\vee}_q$ has a presentation
by generators $z,Z,A$ and relations
\begin{eqnarray*}
&&
\frac{qzA-q^{-1}Az}{q-q^{-1}} = a z^2+a^{-1},
\qquad \qquad
\frac{qZA-q^{-1}AZ}{q-q^{-1}} = a^{-1} Z^2+a,
\\
&&
\qquad
a z^2 -\frac{a^{-1}q-aq^{-1}}{q-q^{-1}}zZ-
\frac{aq-a^{-1}q^{-1}}{q-q^{-1}}Zz + a^{-1}Z^2 = 0.
\end{eqnarray*}
An isomorphism with the presentation in
Theorem
\ref{lem:Uveepres} sends
\begin{eqnarray*}
z &\mapsto& z,
\\
Z &\mapsto& a^2 z + (1-a^2)y^{-1},
\\
A &\mapsto& a^{-1}x+az.
\end{eqnarray*}
The inverse isomorphism sends
\begin{eqnarray*}
x &\mapsto& aA-a^{2}z,
\\
y^{-1} &\mapsto& \frac{az-a^{-1}Z}{a-a^{-1}},
\\
z &\mapsto& z.
\end{eqnarray*}
\end{theorem}
\noindent {\it Proof:} 
Similar to the proof of Theorem
\ref{thm:uveepres3}.
\hfill $\Box$ \\

\noindent We now consider $\nu_x, \nu_y, \nu_z$ and  $\Lambda$ 
from the point of view of $U^\vee_q$.

\begin{lemma}
\label{lem:fail}
The subalgebra $U^\vee_q$ does not contain
any of $\nu_x$,  $\nu_z$, $\Lambda$.
\end{lemma}
\noindent {\it Proof:} 
Write these elements in the basis
for
$U_q({\mathfrak{sl}_2})$ 
 from Lemma \ref{lem:Uqbasis}. 
This can be done using 
(\ref{eq:nux}), (\ref{eq:nuz}) along with
the top left equation in Lemma
\ref{lem:casequit}.
 Evaluate the results using
Corollary
\ref{cor:Uqbasisyi}.
\hfill $\Box$ \\

\noindent In view of Lemma
\ref{lem:fail} we consider the following elements
of 
$U_q({\mathfrak{sl}_2})$:
\begin{eqnarray}
\nu_x y^{-1}, \qquad  \nu_y, \qquad \nu_z y^{-1}, \qquad \Lambda y^{-1}.
\label{eq:four}
\end{eqnarray}

\begin{proposition} 
\label{prop:nxpres1}
Each of the elements
{\rm (\ref{eq:four})} 
is
contained in
$U^\vee_q$.
These elements look as follows in the basis
for $U^\vee_q$
from Corollary
\ref{cor:Uqbasisyi}:
\begin{eqnarray}
&&\nu_x y^{-1}= q^{-1}(y^{-1}-z),
\quad \qquad
\nu_y = q^{-1}(1-xz),
\quad \qquad
\nu_z y^{-1}= q(y^{-1}-x),
\label{eq:nxyietc}
\\
&&\qquad \quad \Lambda y^{-1} = 
q^{-1}(1-xz) + q^{-1}xy^{-1}+q^{-1}y^{-1}z+(q-q^{-1})y^{-2}.
\label{eq:Lamyi}
\end{eqnarray}
\end{proposition}
\noindent {\it Proof:} 
Use (\ref{eq:nux})--(\ref{eq:nuz}) along with
the second equation on the right in 
Lemma
\ref{lem:casequit}.
\hfill $\Box$ \\

\noindent Proposition
\ref{prop:nxpres1} shows how
the elements 
(\ref{eq:four}) 
look in the presentation of $U^\vee_q$ from
Theorem
\ref{lem:Uveepres}. We now show how the elements 
(\ref{eq:four}) 
look
in our other three presentations of $U^\vee_q$.

\begin{proposition} 
\label{prop:nxpres2}
In the presentation of
$U^\vee_q$ from Theorem
\ref{thm:uvpres2},
the elements {\rm (\ref{eq:four})} 
look as follows.
\begin{enumerate}
\item[\rm (i)]
$\nu_x y^{-1}$ is equal to
\begin{eqnarray*}
q^{-1}a^{-2}(y^{-1}-Z).
\end{eqnarray*}
\item[\rm (ii)]
$\nu_y$ is equal to $q^{-1}$ times
\begin{eqnarray*}
1-XZ+(a-a^{-1})a^{-1}y^{-1}Z-(a-a^{-1})aXy^{-1}+(a-a^{-1})^2 y^{-2}.
\end{eqnarray*}
\item[\rm (iii)]
$\nu_z y^{-1}$ is equal to
\begin{eqnarray*}
qa^2(y^{-1}-X).
\end{eqnarray*}
\item[\rm (iv)]
$ \Lambda y^{-1}$ is equal to
\begin{eqnarray*}
q^{-1}(1-XZ) + 
q^{-1}Xy^{-1}
+q^{-1}y^{-1}Z
+(q-q^{-1})y^{-2}.
\end{eqnarray*}
\end{enumerate}
\end{proposition}
\noindent {\it Proof:} (i)--(iii)
By Lemma
\ref{lem:3trans}
we find
$x=a^2X+(1-a^2)y^{-1}$ and
$z=a^{-2}Z+(1-a^{-2})y^{-1}$.
Use these equations to eliminate $x,z$ in
(\ref{eq:nxyietc}).
\\
\noindent (iv) In the equation
(\ref{eq:Lamyi}), apply $\sigma$ to each side.
Use the fact that $\sigma$ fixes each of $\Lambda, y^{-1}$ and sends
$x,z$ to $X,Z$ respectively.
\hfill $\Box$ \\

\begin{proposition} 
\label{prop:nxpres3}
In the presentation of
$U^\vee_q$ from Theorem
\ref{thm:uveepres3},
the elements {\rm (\ref{eq:four})} 
look as follows.
\begin{enumerate}
\item[\rm (i)] $\nu_x y^{-1}$ is equal to $q^{-1}a^{-1}$ times
\begin{eqnarray*}
\frac{a^2X-a^{-2}x}{a-a^{-1}}-A.
\end{eqnarray*}
\item[\rm (ii)] $\nu_y$  is equal to  
\begin{eqnarray*}
a^{-1} \frac{Ax-xA}{q-q^{-1}}.
\end{eqnarray*}
\item[\rm (iii)] $\nu_z y^{-1}$  is equal to 
\begin{eqnarray*}
qa \frac{X-x}{a-a^{-1}}.
\end{eqnarray*}
\item[\rm (iv)] $\Lambda y^{-1}$  is equal to each of
\begin{eqnarray*}
&&
q{ A} \frac{X-x}{a-a^{-1}} +a^{-1}\frac{q-q^{-1}}{a-a^{-1}}x^2
+ \frac{aq^{-1}-a^{-1}q}{a-a^{-1}}xX + q,
\\
&&
q^{-1} \frac{X-x}{a-a^{-1}}{ A} -a^{-1}\frac{q-q^{-1}}{a-a^{-1}}x^2
+ \frac{aq-a^{-1}q^{-1}}{a-a^{-1}}Xx + q^{-1},
\\
&&
q{A} \frac{X-x}{a-a^{-1}} -a\frac{q-q^{-1}}{a-a^{-1}}X^2
+ \frac{aq-a^{-1}q^{-1}}{a-a^{-1}}Xx + q,
\\
&&
q^{-1} \frac{X-x}{a-a^{-1}}{ A} +a\frac{q-q^{-1}}{a-a^{-1}}X^2
+ \frac{aq^{-1}-a^{-1}q}{a-a^{-1}}xX + q^{-1}.
\end{eqnarray*}
\end{enumerate}
\end{proposition}
\noindent {\it Proof:} 
By Lemma
\ref{lem:3trans}
we find $y^{-1}=(aX-a^{-1}x)(a-a^{-1})^{-1}$
and
$z=a^{-1}A-a^{-2}x$. Use these equations to 
eliminate $y^{-1},z$ in
(\ref{eq:nxyietc}),
(\ref{eq:Lamyi}).
Evaluate the results using the relations
in Theorem
\ref{thm:uveepres3}.
\hfill $\Box$ \\

\begin{proposition} 
\label{prop:nxpres4}
In the presentation of
$U^\vee_q$ from Theorem
\ref{thm:uveepres4},
the elements {\rm (\ref{eq:four})} 
look as follows.
\begin{enumerate}
\item[\rm (i)] $\nu_x y^{-1}$ is equal to 
\begin{eqnarray*}
q^{-1}a^{-1} \frac{z-Z}{a-a^{-1}}.
\end{eqnarray*}
\item[\rm (ii)] $\nu_y$  is equal to 
\begin{eqnarray*}
a \frac{zA-Az}{q-q^{-1}}.
\end{eqnarray*}
\item[\rm (iii)] $\nu_z y^{-1}$  is equal to 
$qa$ times
\begin{eqnarray*}
\frac{a^2z-a^{-2}Z}{a-a^{-1}}-A.
\end{eqnarray*}
\item[\rm (iv)] $\Lambda y^{-1}$  is equal to  each of
\begin{eqnarray*}
&&
q^{-1}A \frac{z-Z}{a-a^{-1}} +a\frac{q-q^{-1}}{a-a^{-1}}z^2
+ \frac{aq^{-1}-a^{-1}q}{a-a^{-1}}zZ + q^{-1},
\\
&&
q \frac{z-Z}{a-a^{-1}} A -a\frac{q-q^{-1}}{a-a^{-1}}z^2
+ \frac{aq-a^{-1}q^{-1}}{a-a^{-1}}Zz + q,
\\
&&
q^{-1} A \frac{z-Z}{a-a^{-1}} -a^{-1}\frac{q-q^{-1}}{a-a^{-1}}Z^2
+ \frac{aq-a^{-1}q^{-1}}{a-a^{-1}}Zz + q^{-1},
\\
&&
q \frac{z-Z}{a-a^{-1}} A +a^{-1}\frac{q-q^{-1}}{a-a^{-1}}Z^2
+ \frac{aq^{-1}-a^{-1}q}{a-a^{-1}}zZ + q.
\end{eqnarray*}
\end{enumerate}
\end{proposition}
\noindent {\it Proof:}
By Lemma
\ref{lem:3trans}
we find
$x=aA-a^{2}z$ and
$y^{-1}=(az-a^{-1}Z)(a-a^{-1})^{-1}$.
Use these equations to 
eliminate $x,y^{-1}$ in
(\ref{eq:nxyietc}),
(\ref{eq:Lamyi}),
and evaluate the results using the relations
in Theorem
\ref{thm:uveepres4}. Part (iv) can also be
obtained from
Proposition \ref{prop:nxpres3}(iv), by
applying the antiautomorphism $\dagger$ to
each term in that part of the proposition.
Use the fact that $\dagger$ fixes each of $\Lambda,y^{-1}$
and sends $x,X$ to $Z,z$ respectively.
\hfill $\Box$ \\

\section{Modules for
$U_q({\mathfrak{sl}_2})$ and
$U^{\vee}_q$}

\noindent 
In this section we compare 
the finite-dimensional modules for
$U_q({\mathfrak{sl}_2})$
and $U^{\vee}_q$.
\medskip

\noindent 
Throughout this section $V$ denotes a vector space over
$\F$ with finite positive dimension.
We will use the following fact from linear algebra.
\begin{lemma} 
{\rm \cite[p.~193]{curtis}}.
\label{lem:geninv}
Let $T:V\to V$ denote an invertible linear transformation.
Then there exists a polynomial $f$ that has all coefficients in $\F$
and $T^{-1}=f(T)$.
\end{lemma}

\noindent A module for a given $\F$-algebra is called
{\it decomposable} (resp. {\it semisimple}) whenever it is the direct sum
of two nonzero submodules (resp. direct sum of irreducible submodules).

\medskip
\noindent 
For the moment assume that $V$ is a
$U_q({\mathfrak{sl}_2})$-module. If we restrict
the
$U_q({\mathfrak{sl}_2})$-action
to $U^{\vee}_q$, then  $V$ becomes
a $U^{\vee}_q$-module. 

\begin{lemma} 
\label{lem:uveevsu}
Assume that $V$ is a 
$U_q({\mathfrak{sl}_2})$-module.
For any subspace $W$ of $V$
the following are equivalent:
\begin{enumerate}
\item[\rm (i)]
$W$ is a
$U_q({\mathfrak{sl}_2})$-submodule of $V$;
\item[\rm (ii)]
$W$ is a
$U^\vee_q$-submodule of $V$.
\end{enumerate}
\end{lemma}
\noindent {\it Proof:} 
Assume $W\not=0$; otherwise the result is trivial.
\\
\noindent
${\rm (i)}\Rightarrow {\rm (ii)}$ Since
$U_q({\mathfrak{sl}_2})$ contains
$U^\vee_q$.
\\
\noindent 
${\rm (ii)}\Rightarrow {\rm (i)}$ 
By assumption $W$ is invariant under $y^{-1}$.
Applying Lemma
\ref{lem:geninv}
to $T=y^{-1}$, we see that
$W$ is invariant under $y$. 
\hfill $\Box$ \\


\begin{corollary} 
Assume that $V$ is a $U_q({\mathfrak{sl}_2})$-module.
\begin{enumerate}
\item[\rm (i)] The 
$U_q({\mathfrak{sl}_2})$-module $V$ is irreducible if and only if the
$U^{\vee}_q$-module $V$ is irreducible.
\item[\rm (ii)] The 
$U_q({\mathfrak{sl}_2})$-module $V$ is decomposable if and only if the
$U^{\vee}_q$-module $V$ is decomposable.
\item[\rm (iii)] The 
$U_q({\mathfrak{sl}_2})$-module $V$ is semisimple if and only if
the $U^{\vee}_q$-module $V$ is semisimple.
\end{enumerate}
\end{corollary}
\noindent {\it Proof:} 
Use Lemma
\ref{lem:uveevsu}.
\hfill $\Box$ \\

\begin{note}\label{note:ss}
\rm
Assume that $q$ is not a root of unity,
and that $V$ is a
$U_q({\mathfrak{sl}_2})$-module.
 If $y$ is diagonalizable on $V$
then $V$ is semisimple
\cite[Theorem~2.9]{jantzen}.
Moreover,  
if the characteristic of $\F$ is not 2 
then
$y$ is diagonalizable on $V$
\cite[Proposition~2.3]{jantzen}.
\end{note}

\noindent For the moment assume that $V$ is a
$U^{\vee}_q$-module. We now display some necessary and sufficient
conditions for the 
$U^{\vee}_q$-action on
$V$ 
to extend to a 
$U_q({\mathfrak{sl}_2})$-action on $V$.

\begin{lemma}
\label{lem:ext4}
Assume that $V$ is a 
$U^{\vee}_q$-module. Then the following are equivalent:
\begin{enumerate}
\item[\rm (i)] the action of $U^{\vee}_q$
on $V$ extends to an action of
$U_q({\mathfrak{sl}_2})$ on $V$;
\item[\rm (ii)] the element $y^{-1}$ is invertible on $V$;
\item[\rm (iii)] the element $a^{-1}x-aX$ is invertible on $V$;
\item[\rm (iv)] the element $az-a^{-1}Z$ is invertible on $V$.
\end{enumerate}
\end{lemma}
\noindent {\it Proof:}
${\rm (i)}\Leftrightarrow {\rm (ii)}$ 
By construction.
\\
\noindent ${\rm (ii)}\Leftrightarrow {\rm (iii)}$ 
Use
Lemma
\ref{lem:yi}.
\\
\noindent ${\rm (ii)}\Leftrightarrow {\rm (iv)}$ 
Use
Lemma
\ref{lem:yi}.
\hfill $\Box$ \\

\noindent We comment on linear algebra.
Let $T:V\to V$
denote a linear transformation.
Using $T$ we define the subspaces 
 $V_{inv}$ and $V_{nil}$ of $V$ 
as follows.
For all integers $n\geq 0$
let $V_n$ (resp. $V^{(n)}$)
denote the image (resp. kernel)
of $T^n$ on $V$.
By construction $V_n$ and $V^{(n)}$ are $T$-invariant.
By linear algebra 
${\rm dim}\,V=
{\rm dim}\, V_n + {\rm dim}\, V^{(n)}$.
For
$n\geq 1$ we have
$V_{n} \subseteq  V_{n-1}$ and
$V^{(n-1)} \subseteq  V^{(n)}$.
Define $V_{inv} = \cap_{n=0}^\infty 
V_n$ and
$V_{nil} = \cup_{n=0}^\infty V^{(n)}$.
The subspaces $V_{inv}$ and $V_{nil}$ are $T$-invariant.
By construction
$T$ is invertible 
on
$V_{inv}$ provided $V_{inv}\not=0$.
Moreover $T$ is nilpotent on
$V_{nil}$. By these comments
$V_{inv}\cap
V_{nil}=0$.
Considering the dimensions we obtain
$V=
V_{inv}+
V_{nil}$ (direct sum).
We call $V_{inv}$ (resp.
$V_{nil}$) the {\it invertible part} (resp. {\it nilpotent part})
of $V$ with respect to $T$.

\begin{lemma}
\label{lem:neqone}
Assume that $V$ is a $U^{\vee}_q$-module. Then
the following are
$U^{\vee}_q$-submodules of $V$:
\begin{enumerate}
\item[\rm (i)] the image of $y^{-1}$ on $V$;
\item[\rm (ii)] the kernel of $y^{-1}$ on $V$.
\end{enumerate}
\end{lemma}
\noindent {\it Proof:}  Use
the relations 
(\ref{eq:Uveepres}).
\hfill $\Box$ \\

\begin{lemma}
\label{lem:nil}
Assume that $V$ is a $U^{\vee}_q$-module. Then
for $n \in \N$ the following are
$U^{\vee}_q$-submodules of $V$:
\begin{enumerate}
\item[\rm (i)] the image of $y^{-n}$ on $V$;
\item[\rm (ii)] the kernel of $y^{-n}$ on $V$.
\end{enumerate}
\end{lemma}
\noindent {\it Proof:} (i)
Denote this image  by 
 $V_n$.
Our proof is by induction on $n$.
The case $n=0$ is trivial, so assume 
$n\geq 1$.
By induction $V_{n-1}$ is 
invariant under $U^\vee_q$.
For the $U^{\vee}_q$-module $V_{n-1}$ consider
the image of $y^{-1}$.
By Lemma
\ref{lem:neqone}(i) this image is invariant under
$U^{\vee}_q$.
 By construction this image is 
equal to $V_n$.
Therefore $V_n$ is a 
$U^{\vee}_q$-submodule of $V$.
\\
\noindent (ii)
Denote this kernel  by 
 $V^{(n)}$.
Our proof is by induction on $n$.
The case $n=0$ is trivial, so assume $n\geq 1$.
By induction $V^{(n-1)}$ is 
invariant under $U^\vee_q$.
For the quotient $U^{\vee}_q$-module $V/V^{(n-1)}$
consider the kernel of $y^{-1}$.
By
 Lemma
\ref{lem:neqone}(ii)
this kernel is invariant under
 $U^{\vee}_q$.
By construction  this kernel is equal to
$V^{(n)}/V^{(n-1)}$.
Therefore $V^{(n)}/V^{(n-1)}$ is a 
$U^{\vee}_q$-submodule of
$V/V^{(n-1)}$, and consequently
$V^{(n)}$ is
a $U^{\vee}_q$-submodule of $V$.
\hfill $\Box$ \\

\begin{lemma}
\label{lem:nidec}
Assume that $V$ is a 
$U^\vee_q$-module. Then each of the following
is a 
$U^\vee_q$-submodule of $V$: 
\begin{enumerate}
\item[\rm (i)]
the invertible part 
of $V$ with respect to $y^{-1}$;
\item[\rm (ii)]
the nilpotent part 
of $V$ with respect to $y^{-1}$.
\end{enumerate}
\end{lemma}
\noindent {\it Proof:} Use Lemma
\ref{lem:nil}.
\hfill $\Box$ \\

\noindent We now summarize the situation so far.

\begin{lemma}
\label{lem:sum}
Assume that $V$ is a $U^{\vee}_q$-module.
Let $V_{inv}$ (resp. $V_{nil}$) denote the invertible
part (resp. nilpotent part) of $V$ with respect to $y^{-1}$.
Then 
\begin{eqnarray}
V = V_{inv} + V_{nil} \qquad \qquad 
{\hbox{\rm (direct sum of $U^\vee_q$-modules)}}.
\label{eq:sum}
\end{eqnarray}
Assume that 
 $V_{inv}$ 
 (resp. $V_{nil}$) is nonzero.
 Then on $V_{inv}$
 (resp. $V_{nil}$)
  the action
 of $U^\vee_q$ does (resp. does not) extend to an 
action of $U_q({\mathfrak{sl}_2})$.
\end{lemma}
\noindent {\it Proof:}
Line 
(\ref{eq:sum}) is from the construction and
Lemma
\ref{lem:nidec}.
To get the remaining assertions,
apply parts (i), (ii) of Lemma
\ref{lem:ext4} to
 $V_{inv}$ and
$V_{nil}$.
\hfill $\Box$ \\

\noindent A module for a given $\F$-algebra is called
{\it indecomposable} whenever it is not decomposable.

\begin{corollary}
\label{cor:nilorinv}
Assume that $V$ is an
indecomposable $U^\vee_q$-module.
Then on $V$ the element $y^{-1}$ is either nilpotent or invertible.
\end{corollary}
\noindent {\it Proof:}
Referring to (\ref{eq:sum}),
either 
$V_{inv}=0$ or
$V_{nil}=0$.
\hfill $\Box$ \\

\noindent For the moment assume that $V$ is a
$U_q({\mathfrak{sl}_2})$-module.
Then on the 
$U^\vee_q$-module $V$ the element $y^{-1}$ is invertible.
We now give some  
examples of a 
$U^\vee_q$-module on which
$y^{-1}$ is nilpotent.

\begin{definition}\rm
Let $\mathbb A$ denote the $\F$-algebra defined by generators
$u$, $v$ and one relation
\begin{eqnarray*}
\frac{quv-q^{-1}vu}{q-q^{-1}}= 1.
\end{eqnarray*}
\end{definition}

\begin{lemma}
\label{lem:utoa}
There exists an $\F$-algebra homomorphism
$U^{\vee}_q \to \mathbb A$ that sends
\begin{eqnarray*}
x \mapsto v, \qquad \qquad y^{-1} \mapsto 0, \qquad \qquad
z \mapsto u.
\end{eqnarray*}
\end{lemma}
\noindent {\it Proof:}
Use Theorem
\ref{lem:Uveepres}.
\hfill $\Box$ \\

\begin{example}\rm 
Assume that $V$ is an
$\mathbb A$-module.
If we pull back the $\mathbb A$-module structure
via the homomorphism $U^{\vee}_q \to \mathbb A$
from Lemma
\ref{lem:utoa}, then $V$ becomes a $U^{\vee}_q$-module on which
$y^{-1}$ is zero. 
\end{example}

\section{$U^\vee_q$-modules from tridiagonal pairs}
In this section we use tridiagonal pairs
of $q$-Racah type to construct examples of
finite-dimensional $U^\vee_q$-modules.


\medskip

\noindent We now recall the notion of a tridiagonal pair
\cite{TD00}.
For background information on this topic, we refer the reader to
\cite{TD00, IT:qRacah,
TDclass,
NT:muqrac, madrid}.

\begin{definition}
\label{def:tdp}
\rm
\cite[Definition~1.1]{TD00}.
Let $V$ denote a vector space over $\F$ with finite positive
dimension. By a {\it tridiagonal pair on $V$}, we mean an ordered pair
of $\F$-linear maps ${\bf A}:V\to V$ and
 ${\bf A}^*:V\to V$ that satisfy the following conditions.
\begin{enumerate}
\item[\rm (i)] Each of ${\bf A}$, ${\bf A^*}$ is diagonalizable on $V$.
\item[\rm (ii)] There exists an ordering
$\lbrace V_i\rbrace_{i=0}^d$ of the eigenspaces  of $\bf A$ such that
\begin{eqnarray}
{\bf A^*}V_i \subseteq V_{i-1}+V_i+V_{i+1}
\qquad \qquad (0 \leq i \leq d),
\label{eq:TD1}
\end{eqnarray}
where $V_{-1}=0$ and $V_{d+1}=0$.
\item[\rm (iii)] There exists an ordering
$\lbrace V^*_i\rbrace_{i=0}^\delta$ of the eigenspaces  of $\bf A^*$ such that
\begin{eqnarray}
{\bf A}V^*_i \subseteq V^*_{i-1}+V^*_i+V^*_{i+1}
\qquad \qquad (0 \leq i \leq \delta),
\label{eq:TD2}
\end{eqnarray}
where $V^*_{-1}=0$ and $V^*_{\delta+1}=0$.
\item[\rm (iv)] There does not exist a subspace $W \subseteq V$ such that
${\bf A}W\subseteq W$,
${\bf A^*}W\subseteq W$,
$W\not=0$,
$W\not=V$.
\end{enumerate}
\end{definition}

\begin{note}\rm According to a common notational convention ${\bf A^*}$
denotes the conjugate-transpose of $\bf A$. We are not using this
convention. For a tridiagonal pair $\bf A, \bf A^*$ the $\F$-linear
maps $\bf A$ and $\bf A^*$ are arbitrary subject to (i)--(iv) above.
\end{note}

\medskip
\noindent In order to motivate our results, 
we summarize some facts about tridiagonal pairs. For
the rest of this section, fix a tridiagonal pair
$\bf A, \bf A^*$ on $V$ as in Definition
\ref{def:tdp}.
By \cite[Lemma~4.5]{TD00} the integers
$d$ and $\delta$ from (ii) and (iii) are equal;
we call this common value the {\it diameter} of the pair.
To avoid trivialities, we always assume that $d\geq 1$.
An ordering of the eigenspaces of $\bf A$ (resp. $\bf A^*$)
is said to be {\em standard} whenever it satisfies
(\ref{eq:TD1})
 (resp. (\ref{eq:TD2})).
 We comment on the uniqueness of the standard ordering.
 Let $\lbrace V_i\rbrace_{i=0}^d$ 
 denote a standard ordering of the eigenspaces of $\bf A$.
 By \cite[Lemma~2.4]{TD00},
  the ordering $\lbrace V_{d-i} \rbrace_{i=0}^d$ is also standard and no further
   ordering
   is standard.
   A similar result holds for the eigenspaces of $\bf A^*$.
   For the rest of this section fix a standard ordering
   $\lbrace V_i \rbrace_{i=0}^d$ (resp.
   $ \lbrace V^*_i \rbrace_{i=0}^d$)
   of the eigenspaces
    of $ \bf A$ (resp. $\bf A^*$).
    For $0 \leq i \leq d$ let
    $\theta_i$
    (resp. $\theta^*_i$)
    denote the eigenvalue of
    $\bf A$
    (resp.  $ \bf A^*$) associated with
    $V_i$
    (resp. $V^*_i$).
   By construction $\lbrace \theta_i\rbrace_{i=0}^d$ are mutually
   distinct scalars in $\F$, and
   $\lbrace \theta^*_i\rbrace_{i=0}^d$ are mutually
   distinct scalars in $\F$.
    By
     \cite[Theorem~11.1]{TD00}
      the expressions
      \begin{eqnarray*}
      \frac{\theta_{i-2}-\theta_{i+1}}{\theta_{i-1}-\theta_i},  \qquad\qquad
        \frac{\theta^*_{i-2}-\theta^*_{i+1}}{\theta^*_{i-1}-\theta^*_i}
	\end{eqnarray*}
	are equal and independent of $i$ for $2 \leq i \leq d-1$.
\medskip

\noindent We recall the split decomposition \cite[Section~4]{TD00}.
For $0 \leq i \leq d$ define
\begin{eqnarray*}
U_i = (V^*_0+V^*_1+\cdots + V^*_i)\cap
(V_i+V_{i+1}+\cdots + V_d).
\end{eqnarray*}
For notational convenience define $U_{-1}=0$ and
$U_{d+1}=0$.
By \cite[Theorem~4.6]{TD00} we have
$V = \sum_{i=0}^d U_i$ (direct sum). 
Moreover for $0 \leq i \leq d$,
\begin{eqnarray*}
({\bf A}-\theta_i I)U_i\subseteq U_{i+1},
\qquad \qquad 
({\bf A^*}-\theta^*_i I)U_i\subseteq U_{i-1}.
\end{eqnarray*}
For $0 \leq i \leq d$ define
\begin{eqnarray*}
U^\Downarrow_i = (V^*_0+V^*_1+\cdots + V^*_i)\cap
(V_0+V_1+\cdots + V_{d-i}).
\end{eqnarray*}
For notational convenience define $U^\Downarrow_{-1}=0$ and
$U^\Downarrow_{d+1}=0$.
By \cite[Theorem~4.6]{TD00} we have
$V = \sum_{i=0}^d U^\Downarrow_i$ (direct sum). 
Moreover  for $0 \leq i \leq d$,
\begin{eqnarray*}
({\bf A}-\theta_{d-i} I)U^\Downarrow_i\subseteq U^\Downarrow_{i+1},
\qquad \qquad 
({\bf A^*}-\theta^*_i I)U^\Downarrow_i\subseteq U^\Downarrow_{i-1}.
\end{eqnarray*}
\noindent We recall the definition of $q$-Racah type. 
Following
\cite[Definition~3.1]{IT:qRacah} and
\cite[Definition~5.1]{hwh},
we say that
$\bf A, \bf A^*$ has {\it $q$-Racah type} whenever
there exist nonzero $a,b \in \F$ such that  both
\begin{eqnarray*}
\theta_i = aq^{d-2i} + a^{-1}q^{2i-d},
\qquad \qquad 
\theta^*_i = bq^{d-2i} + b^{-1}q^{2i-d}
\end{eqnarray*}
for $0 \leq i \leq d$.
For the rest of this section assume that
$\bf A, \bf A^*$ has $q$-Racah type.
For $1 \leq i \leq d$ we have $q^{2i}\not=1$; otherwise
$\theta_0=\theta_i$.
Also $a^2\not=1$; otherwise $\theta_0=\theta_d$.
Similarly $b^2\not=1$.
We now recall the maps $K,B$.

\begin{definition}\rm
\cite[Definitions~3.1, 3.2]{bockting2}.
Define an $\F$-linear map $K:V\to V$ such that for
$0 \leq i \leq d$, $U_i$ is an eigenspace of $K$ with
eigenvalue $q^{d-2i}$. Thus
\begin{eqnarray*}
(K-q^{d-2i}I)U_i=0 \qquad \qquad (0 \leq i \leq d).
\end{eqnarray*}
Define an $\F$-linear map $B:V\to V$ such that for
$0 \leq i \leq d$, $U^\Downarrow_i$ is an eigenspace of $B$ with
eigenvalue $q^{d-2i}$. Thus
\begin{eqnarray*}
(B-q^{d-2i}I)U^\Downarrow_i=0 \qquad \qquad (0 \leq i \leq d).
\end{eqnarray*}
\end{definition}

\begin{lemma} 
{\rm \cite[Lemma~3.6, Theorem~9.9]{bockting2}}.
\label{lem:bockting} 
The maps $K,B,\bf A$ satisfy
\begin{eqnarray*}
&&
\frac{qK {\bf A} - q^{-1} {\bf A} K}{q-q^{-1}} = a K^2 + a^{-1} I,
\qquad \qquad 
\frac{qB {\bf A} - q^{-1} {\bf A} B}{q-q^{-1}} = a^{-1} B^2 + a I,
\\
&&
\qquad \quad
aK^2 - \frac{a^{-1} q-aq^{-1}}{q-q^{-1}} KB 
- \frac{a q-a^{-1} q^{-1}}{q-q^{-1}} BK  
+
a^{-1} B^2 = 0. 
\end{eqnarray*}
\end{lemma}

\noindent We are done with our summary.
We now use ${\bf A, \bf A^*}$ to construct two
$U^\vee_q$-module structures on $V$.
Here is the first one.

\begin{theorem}
\label{thm:1}
There exists a unique
$U^\vee_q$-module structure on $V$
for which the generators $z,Z,A$ act as follows:

\begin{center}
\begin{tabular}{c|c c c}
{\rm generator} & $z$ &  $Z$ & $A$
\\
\hline
{\rm action on $V$} & $K$ & $B$ &$\bf A$ 
\end{tabular}
 \end{center}

\noindent On the
$U^\vee_q$-module $V$ 
the elements $x, X, y^{-1}$ act as follows:

\begin{center}
\begin{tabular}{c|c c c}
{\rm element} & $x$ &  $X$ & $y^{-1}$
\\
\hline
{\rm action on $V$} & $a {\bf A}-a^2K$ & $a^{-1}{\bf A}-a^{-2}B$ &
$\frac{aK-a^{-1}B}{a-a^{-1}}$ 
\end{tabular}
 \end{center}

\end{theorem}
\noindent {\it Proof:}
To see that the given $U^\vee_q$-module structure exists,
compare the three equations in
Lemma
\ref{lem:bockting} 
with the defining relations for
$U^\vee_q$ given in 
Theorem
\ref{thm:uveepres4}. The 
$U^\vee_q$-module structure is unique since
$z,Z,A$ generate $U^\vee_q$.
The actions of $x$ and $y^{-1}$ on $V$ are obtained
using the identifications involving those elements
given in Theorem
\ref{thm:uveepres4}. The action of $X$ on $V$ is
found using $A=aX+a^{-1}Z$.
\hfill $\Box$ \\

\begin{proposition}
\label{thm:yiext}
On the 
$U^\vee_q$-module $V$ 
from Theorem
\ref{thm:1},
the elements
$\nu_x y^{-1}$, $\nu_y$, $\nu_z y^{-1}$,
$\Lambda y^{-1}$ act as follows.
\begin{enumerate}
\item[\rm (i)]
$\nu_x y^{-1}$ acts as
\begin{eqnarray*}
q^{-1}a^{-1}\frac{K-B}{a-a^{-1}}.
\end{eqnarray*}
\item[\rm (ii)] 
$\nu_y$ acts as 
\begin{eqnarray*}
a\,\frac{K{\bf A}-{\bf A}K}{q-q^{-1}}.
\end{eqnarray*}
\item[\rm (iii)] 
$\nu_z y^{-1}$ acts as $qa$ times
\begin{eqnarray*}
\frac{a^2K-a^{-2}B}{a-a^{-1}}-{\bf A}.
\end{eqnarray*}
\item[\rm (iv)] 
$\Lambda y^{-1}$ acts as  each of
\begin{eqnarray*}
&&
q^{-1}{\bf A} \frac{K-B}{a-a^{-1}} +a\frac{q-q^{-1}}{a-a^{-1}}K^2
+ \frac{aq^{-1}-a^{-1}q}{a-a^{-1}}KB + q^{-1}I,
\\
&&
q \frac{K-B}{a-a^{-1}}{\bf A} -a\frac{q-q^{-1}}{a-a^{-1}}K^2
+ \frac{aq-a^{-1}q^{-1}}{a-a^{-1}}BK + qI,
\\
&&
q^{-1}{\bf A} \frac{K-B}{a-a^{-1}} -a^{-1}\frac{q-q^{-1}}{a-a^{-1}}B^2
+ \frac{aq-a^{-1}q^{-1}}{a-a^{-1}}BK + q^{-1}I,
\\
&&
q \frac{K-B}{a-a^{-1}}{\bf A} +a^{-1}\frac{q-q^{-1}}{a-a^{-1}}B^2
+ \frac{aq^{-1}-a^{-1}q}{a-a^{-1}}KB + qI.
\end{eqnarray*}
\end{enumerate}
\end{proposition}
\noindent {\it Proof:}
Apply Proposition
\ref{prop:nxpres4}
with $z=K$, $Z=B$, $A={\bf A}$.
\hfill $\Box$ \\

\noindent We have a comment.
\begin{lemma}
\label{lem:expinv}
{\rm 
 \cite[Lemma~9.7]{bockting2}}.
The map $aI-a^{-1}BK^{-1}$ is invertible.
\end{lemma}

\begin{lemma} 
\label{lem:extendable1}
On the 
$U^\vee_q$-module $V$ from Theorem
\ref{thm:1}, the element
$y^{-1}$ is invertible.
\end{lemma}
\noindent {\it Proof:}
By Theorem
\ref{thm:1}, $y^{-1}$ acts on $V$
as a nonzero scalar multiple of
 $aK-a^{-1}B$.
 We have $aK-a^{-1}B=(aI-a^{-1}BK^{-1})K$,
and  
$aI-a^{-1}BK^{-1}$ is invertible by Lemma
\ref{lem:expinv}. 
The result follows.
\hfill $\Box$ \\

\noindent We now bring in the Bockting operator $\Psi$
associated with $\bf A,\bf A^*$.
This was introduced in
\cite{twocom} and investigated further in
\cite{bockting2}. Following
\cite{bockting2}
we will work with the normalized version
\begin{eqnarray*}
\psi = (q-q^{-1})(q^d-q^{-d})\Psi.
\end{eqnarray*}
One feature of $\psi$ is 
that $\psi U_i \subseteq U_{i-1}$ and
$\psi U^\Downarrow_i \subseteq U^\Downarrow_{i-1}$
for $0 \leq i \leq d$ \cite[Lemma~11.2,~Corollary~15.3]{twocom}.
By \cite[Theorem~9.8]{bockting2},
\begin{eqnarray}
\label{eq:step1}
\psi &=& q^{-1}(I-BK^{-1})(aI-a^{-1}BK^{-1})^{-1}.
\end{eqnarray}

\begin{proposition}
\label{thm:psiext}
The 
$U^\vee_q$-module $V$ from Theorem
\ref{thm:1} extends to a 
$U_q({\mathfrak{sl}_2})$-module.
On the
$U_q({\mathfrak{sl}_2})$-module $V$ we have
 $\nu_x=a^{-1}\psi$.
\end{proposition}
\noindent {\it Proof:}
The first assertion follows from
Lemma \ref{lem:ext4}(i),(ii)
and
Lemma \ref{lem:extendable1}.
We now prove the second assertion.
Using Theorem
\ref{thm:1} the following holds on $V$:
\begin{eqnarray}
\label{eq:step2}
y^{-1} = \frac{aK-a^{-1}B}{a-a^{-1}} 
=
\frac{aI-a^{-1}BK^{-1}}{a-a^{-1}}K. 
\end{eqnarray}
Using
(\ref{eq:step1}), 
(\ref{eq:step2}) we find
$\psi y^{-1} = q^{-1}(K-B)(a-a^{-1})^{-1}$, which is equal to
$a\nu_x y^{-1}$
by Proposition
\ref{thm:yiext}(i).
Now invoking
Lemma \ref{lem:extendable1}
we see that $\nu_x  = a^{-1}\psi$ on $V$.
\hfill $\Box$ \\

\begin{proposition}
\label{prop:yacts}
On the
$U_q({\mathfrak{sl}_2})$-module $V$ from
Proposition
\ref{thm:psiext},
the action of $y$ coincides with each of
the following:
\begin{eqnarray*}
K^{-1}(1-a^{-1}q\psi),
\qquad \qquad 
(1-a^{-1}q^{-1}\psi)K^{-1},
\\
B^{-1}(1-aq\psi),
\qquad \qquad 
(1-aq^{-1}\psi)B^{-1}.
\end{eqnarray*}
\end{proposition}
\noindent {\it Proof:}
The four displayed expressions are equal by
\cite[Proposition~9.2]{bockting2}.
We show that $y=K^{-1}(1-a^{-1}q\psi)$ on $V$.
Recall that
$\nu_x = q^{-1}(1-zy)$. 
We have seen that on $V$,
$\nu_x =a^{-1}\psi$
and $z=K$. By these comments
$a^{-1}\psi =
q^{-1}(1-Ky)$ on $V$.
Solve this equation for $y$ to find that
$y=K^{-1}(1-a^{-1}q\psi)$ on $V$.
\hfill $\Box$ \\


\begin{proposition} 
\label{prop:Lamacts}
On 
the 
$U_q({\mathfrak{sl}_2})$-module $V$ from
Proposition
\ref{thm:psiext}, the action of $\Lambda$ 
coincides with each of the following:
\begin{eqnarray}
&&({\bf A}-aK-a^{-1}K^{-1})\psi+qK+q^{-1}K^{-1},
\label{eq:c1}
\\
&&\psi({\bf A}-aK-a^{-1}K^{-1})+q^{-1}K+qK^{-1},
\label{eq:c2}
\\
&&({\bf A}-a^{-1}B-a B^{-1})\psi+qB+q^{-1}B^{-1},
\label{eq:c3}
\\
&&\psi({\bf A}-a^{-1}B-aB^{-1})+q^{-1}B+qB^{-1}.
\label{eq:c4}
\end{eqnarray}
\end{proposition}
\noindent {\it Proof:}
The expressions 
(\ref{eq:c1})--(\ref{eq:c4})
are equal by
\cite[Lemma~9.1]{bockting2}.
Using in order
 Proposition
\ref{prop:yacts},
Theorem
\ref{thm:1},
Proposition \ref{thm:psiext},
and line 
(\ref{eq:nux}), we find that on
 $V$ the 
expression (\ref{eq:c1}) is equal to
\begin{eqnarray*}
&&({\bf A}-aK)\psi + qK + q^{-1}K^{-1}(1-a^{-1}q\psi)
\\
&& =
({\bf A}-aK)\psi + qK + q^{-1}y 
\\
&&= x \nu_x + q z + q^{-1} y
\\
&&= qx + q^{-1}y + qz-qxyz,
\end{eqnarray*}
which is equal to $\Lambda$ by
Lemma
\ref{lem:casequit}.
The result follows.
\hfill $\Box$ \\

\noindent We are done discussing the $U^\vee_q$-module
$V$ from Theorem
\ref{thm:1}. We now consider another $U^\vee_q$-module structure on $V$.

\begin{lemma} 
\label{lem:bocktingalt} 
The maps $K,B,\bf A$ satisfy
\begin{eqnarray*}
\frac{q {\bf A}K^{-1} - q^{-1} K^{-1}{\bf A} }{q-q^{-1}} = a^{-1} K^{-2}+aI,
\qquad \quad 
\frac{q {\bf A}B^{-1} - q^{-1} B^{-1}{\bf A}}{q-q^{-1}} = a B^{-2}+a^{-1} I.
\end{eqnarray*}
\end{lemma}
\noindent {\it Proof:}
These are reformulations of the first two equations in
Lemma \ref{lem:bockting}.
\hfill $\Box$ \\

\begin{lemma}
\label{lem:bocktingalt2} 
{\rm \cite[Theorem~9.10]{bockting2}.}
The maps $K,B$ satisfy
\begin{eqnarray*}
a^{-1}K^{-2} - \frac{a^{-1} q-aq^{-1}}{q-q^{-1}} K^{-1}B^{-1} 
- \frac{a q-a^{-1} q^{-1}}{q-q^{-1}} B^{-1}K^{-1}  
+
a B^{-2} = 0. 
\end{eqnarray*}
\end{lemma}

\begin{theorem}
\label{thm:2}
There exists a 
$U^\vee_q$-module structure on $V$
for which the generators $x,X,A$ act as follows:

\begin{center}
\begin{tabular}{c|c c c}
{\rm generator} & $x$ &  $X$ & $A$
\\
\hline
{\rm action on $V$} & $K^{-1}$ & $B^{-1}$ &$\bf A$ 
\end{tabular}
 \end{center}

\noindent On this
$U^\vee_q$-module $V$ 
the elements $z, Z, y^{-1}$ act as follows:

\begin{center}
\begin{tabular}{c|c c c}
{\rm element} & $z$ &  $Z$ & $y^{-1}$
\\
\hline
{\rm action on $V$} & $a^{-1} {\bf A}-a^{-2}K^{-1}$ &
$a{\bf A}-a^2B^{-1}$ &
$\frac{a^{-1}K^{-1}-aB^{-1}}{a^{-1}-a}$ 
\end{tabular}
 \end{center}

\end{theorem}
\noindent {\it Proof:}
To see that the given $U^\vee_q$-module structure exists,
compare the equations in
Lemmas \ref{lem:bocktingalt},
\ref{lem:bocktingalt2} 
with the defining relations for
$U^\vee_q$ given in 
Theorem
\ref{thm:uveepres3}. The 
$U^\vee_q$-module structure is unique since
$x,X,A$ generate $U^\vee_q$.
The actions of $y^{-1}$ and $z$ on $V$ are obtained
using the identifications involving those elements
given in Theorem
\ref{thm:uveepres3}. The action of $Z$ on $V$ is
found using $A=aX+a^{-1}Z$.
\hfill $\Box$ \\

\begin{proposition} 
\label{prop:nxpresBK}
On the
$U^\vee_q$-module $V$ from Theorem
\ref{thm:2},
the elements 
$\nu_x y^{-1}$, $\nu_y$, $\nu_z y^{-1}$,
$\Lambda y^{-1}$
act as follows.
\begin{enumerate}
\item[\rm (i)] $\nu_x y^{-1}$ acts as $q^{-1}a^{-1}$ times
\begin{eqnarray*}
\frac{a^2B^{-1}-a^{-2}K^{-1}}{a-a^{-1}}-{\bf A}.
\end{eqnarray*}
\item[\rm (ii)] $\nu_y$  acts as
\begin{eqnarray*}
a^{-1} \frac{{\bf A} K^{-1}-K^{-1}{\bf A}}{q-q^{-1}}.
\end{eqnarray*}
\item[\rm (iii)] $\nu_z y^{-1}$  acts as
\begin{eqnarray*}
qa \frac{B^{-1}-K^{-1}}{a-a^{-1}}.
\end{eqnarray*}
\item[\rm (iv)] $\Lambda y^{-1}$  acts as each of
\begin{eqnarray*}
&&
q{ \bf A} \frac{B^{-1}-K^{-1}}{a-a^{-1}} +a^{-1}\frac{q-q^{-1}}{a-a^{-1}}K^{-2}
+ \frac{aq^{-1}-a^{-1}q}{a-a^{-1}}K^{-1}B^{-1} + qI,
\\
&&
q^{-1} \frac{B^{-1}-K^{-1}}{a-a^{-1}}{\bf A} -a^{-1}\frac{q-q^{-1}}{a-a^{-1}}K^{-2}
+ \frac{aq-a^{-1}q^{-1}}{a-a^{-1}}B^{-1}K^{-1} + q^{-1}I,
\\
&&
q{\bf A} \frac{B^{-1}-K^{-1}}{a-a^{-1}} -a\frac{q-q^{-1}}{a-a^{-1}}B^{-2}
+ \frac{aq-a^{-1}q^{-1}}{a-a^{-1}}B^{-1}K^{-1} + qI,
\\
&&
q^{-1} \frac{B^{-1}-K^{-1}}{a-a^{-1}}{ \bf A} +a\frac{q-q^{-1}}{a-a^{-1}}B^{-2}
+ \frac{aq^{-1}-a^{-1}q}{a-a^{-1}}K^{-1}B^{-1} + q^{-1}I.
\end{eqnarray*}
\end{enumerate}
\end{proposition}
\noindent {\it Proof:}
Apply Proposition \ref{prop:nxpres3}
with $x=K^{-1}$,
$X=B^{-1}$,
$A={\bf A}$.
\hfill $\Box$ \\

\begin{lemma} On the 
$U^\vee_q$-module $V$ from Theorem
\ref{thm:2}, the element
$y^{-1}$ is invertible.
\end{lemma}
\noindent {\it Proof:}
Similar to the proof of
Lemma
\ref{lem:extendable1}.
\hfill $\Box$ \\

\begin{proposition} 
\label{thm:uq2ext}
The 
$U^\vee_q$-module $V$ from Theorem
\ref{thm:2} extends to a 
$U_q({\mathfrak{sl}_2})$-module.
On the 
$U_q({\mathfrak{sl}_2})$-module $V$ we have
$\nu_z=a\psi$.
\end{proposition}
\noindent {\it Proof:}
It suffices to show that $y^{-1} \nu_z = a y^{-1} \psi$ on $V$.
To show this, use $y^{-1}\nu_z = q^{-2}\nu_z y^{-1}$ 
and an argument similar to 
 the proof of
Proposition \ref{thm:psiext}.
\hfill $\Box$ \\

\begin{proposition}
On the
$U_q({\mathfrak{sl}_2})$-module $V$
from Proposition
\ref{thm:uq2ext}, the action of $y$ coincides with
each of the following:
\begin{eqnarray*}
K(1-aq^{-1}\psi),
\qquad \qquad 
(1-aq\psi)K,
\\
B(1-a^{-1}q^{-1}\psi),
\qquad \qquad 
(1-a^{-1}q\psi)B.
\end{eqnarray*}
\end{proposition}
\noindent {\it Proof:}
The four displayed expressions are equal by
\cite[Proposition~9.2]{bockting2}.
To show that $y=K(1-aq^{-1}\psi)$ on $V$,
proceed as in the proof of
Proposition \ref{prop:yacts}.
\hfill $\Box$ \\

\begin{proposition}
\label{prop:Lamacts2}
On the $U_q({\mathfrak{sl}_2})$-module $V$
from Proposition
\ref{thm:uq2ext}, the action of $\Lambda $ coincides with
each of
 {\rm (\ref{eq:c1})--(\ref{eq:c4})}.
\end{proposition}
\noindent {\it Proof:}
Similar to the proof of
Proposition
\ref{prop:Lamacts}.
\hfill $\Box$ \\

\begin{corollary} The following coincide:
\begin{enumerate} 
\item[\rm (i)]
the action of $\Lambda$ on the
$U_q({\mathfrak{sl}_2})$-module $V$ 
from Proposition
\ref{thm:psiext};
\item[\rm (ii)]
the action of $\Lambda$ on the
$U_q({\mathfrak{sl}_2})$-module $V$ 
from Proposition
\ref{thm:uq2ext}.
\end{enumerate}
\end{corollary}
\noindent {\it Proof:}
Compare Proposition
\ref{prop:Lamacts} and
Proposition
\ref{prop:Lamacts2}.
\hfill $\Box$ \\

\noindent To summarize, in each of Theorems
\ref{thm:1} and
\ref{thm:2} we 
displayed a $U^\vee_q$-module
structure on $V$. Using Propositions
\ref{thm:psiext} and
\ref{thm:uq2ext},
we extended each 
$U^\vee_q$-module $V$ to a 
$U_q({\mathfrak{sl}_2})$-module $V$.
Holding on these $U_q({\mathfrak{sl}_2})$-modules
are all the
$U_q({\mathfrak{sl}_2})$ relations
from Sections 2--9, 
such as in
Corollary \ref{cor:yAy},
Lemma \ref{lem:CasAy}, Proposition
\ref{prop:tdrel}, and Lemma
\ref{lem:yAy}.
We expect that this information
will be useful in future developments concerning
the theory of tridiagonal pairs.

\section{Directions for future research}

\noindent In this section we make some suggestions for
future research.
\medskip

\noindent The following problem is motivated by
Lemma \ref{lem:qsymrel}.
Note that the relations
displayed in that lemma
are invariant under the
move $q \mapsto q^{-1}$.

\begin{problem}\rm Consider the $\F$-algebra
$\mathbb U$
defined by generators $x,y,z$ and the relations
from Lemma
\ref{lem:qsymrel}.
Find a basis for the $\F$-vector space $\mathbb U$. 
Note that
the symmetric group $S_3$ acts on 
$\mathbb U$
as a group of automorphisms, by permuting
$x,y,z$. Observe that there exists an
$\F$-algebra homomorphism
$\mathbb U \to 
U_q({\mathfrak{sl}_2})$ that sends
$x\mapsto x$,
$y\mapsto y$,
$z\mapsto z$; let $J_q$ denote the kernel.
Similarly 
there exists an
$\F$-algebra homomorphism
$\mathbb U \to 
U_{q^{-1}}({\mathfrak{sl}_2})$ that sends
$x\mapsto x$,
$y\mapsto y$,
$z\mapsto z$; let $J_{q^{-1}}$ denote the kernel.
Describe the ideals $J_q\cap J_{q^{-1}}$ and
$J_q +J_{q^{-1}}$.
\end{problem}

\begin{problem}\rm 
Find all the $\Z_2$-gradings of
$U_q({\mathfrak{sl}_2})$.
\end{problem}

\begin{problem}\rm 
By Proposition
\ref{prop:gensets}(ii)
the $\F$-algebra $U'_{even}$ is generated by
 $\nu_x,\nu_y,\nu_z$.
Using these generators, 
find a presentation of $U'_{even}$
by generators and relations.
Investigate $U'_{odd}$ as a $U'_{even}$-module.
\end{problem}

\noindent To motivate the next problem we have some comments.
Observe that for $\varepsilon \in \lbrace 1,-1\rbrace$
there exists an $\F$-algebra homomorphism
$U_q({\mathfrak{sl}_2})\to \F$ that sends
$x\mapsto \varepsilon $, $y\mapsto \varepsilon $, $z\mapsto \varepsilon $.
For these two homomorphisms
the intersection of their kernels will
be denoted by $J$. 
Note that $J$ is a  2-sided
ideal of 
$U_q({\mathfrak{sl}_2})$.

\begin{problem}
\label{prob:J}
\rm 
Show that the above 2-sided ideal $J$ is generated by
 $\nu_x,\nu_y,\nu_z$. 
Also show that
\begin{eqnarray*}
U_q({\mathfrak{sl}_2}) = J+ \F1 + \F y 
\qquad \qquad
{\mbox{\rm (direct sum of vector spaces).}}
\end{eqnarray*}
\end{problem}

\begin{problem} \rm 
Show that the center of
$U^\vee_q$ is equal to $\F 1$, provided that
$q$ is not a root of unity.
\end{problem}

\noindent 
The following problem is motivated by 
Corollary
\ref{cor:nilorinv}.
\begin{problem}\rm Classify up to isomorphism the
finite-dimensional
indecomposable $U^{\vee}_q$-modules on which
$y^{-1}$ is nilpotent.
\end{problem}


\begin{note}\rm
We mention some infinite-dimensional
$U^\vee_q$-modules on which $y^{-1}$ is nilpotent.
Recall the basis 
(\ref{eq:xyzbasis}) 
for
$U_q({\mathfrak{sl}_2})$.
For $i \in \Z$ let $U_i$ denote the
subspace of
$U_q({\mathfrak{sl}_2})$  spanned by those
basis elements  $x^ry^s z^t$ such that $s \leq i$.
Thus $U_0 = U^\vee_q$. 
We have
$U_{i-1}\subseteq U_i$ for $i \in \Z$. Moreover
$U_q({\mathfrak{sl}_2}) = \cup_{i \in \Z} U_i$.
One checks that 
\begin{eqnarray*}
x U_i \subseteq U_i, \qquad
z U_i \subseteq U_i, \qquad
y U_i = U_{i+1}, \qquad 
y^{-1} U_i = U_{i-1}, \qquad i \in \Z.
\end{eqnarray*}
Moreover $U_i U_j \subseteq U_{i+j}$ for $i,j\in \Z$.
Thus the sequence $\lbrace U_i\rbrace_{i \in \Z}$ is
a $\Z$-filtration of 
$U_q({\mathfrak{sl}_2})$.
View 
$U_q({\mathfrak{sl}_2})$ as a 
$U^\vee_q$-module 
such that $\xi$ sends $v\mapsto \xi v$ for all
$\xi \in U^\vee_q$ and
$v\in U_q({\mathfrak{sl}_2})$.
For $i \in \Z$ the subspace
$U_i$ 
is a $U^\vee_q$-submodule.
For 
 $i,j\in \Z$ such that $j<i$, consider  the quotient
$U^\vee_q$-module
$U_i /U_j$. For $0 \leq r \leq i-j$
the image (resp. kernel) of $y^{-r}$ is
$U_{i-r} /U_j$ (resp. 
$U_{j+r} /U_j$).
In particular $y^{j-i}$ is zero on
$U_i/U_j$, so
$y^{-1}$ is nilpotent on $U_i/U_j$.
\end{note}

\section{Acknowledgement} 
The authors thank Kazumasa Nomura and
Tatsuro Ito for helpful discussions.

\bigskip

\noindent Sarah Bockting-Conrad \hfil\break
\noindent Department of Mathematics \hfil\break
\noindent University of Wisconsin \hfil\break
\noindent 480 Lincoln Drive \hfil\break
\noindent Madison, WI 53706-1388 USA \hfil\break
\noindent email: {\tt bockting@math.wisc.edu}\hfil\break 

\bigskip

\noindent Paul Terwilliger \hfil\break
\noindent Department of Mathematics \hfil\break
\noindent University of Wisconsin \hfil\break
\noindent 480 Lincoln Drive \hfil\break
\noindent Madison, WI 53706-1388 USA \hfil\break
\noindent email: {\tt terwilli@math.wisc.edu }\hfil\break

\end{document}